\documentclass[12pt]{article}
\usepackage{amssymb,epsfig,times,graphicx}
\usepackage{amsmath}

\newcommand{\sign}{\mbox{sign}\hskip 0.5truemm}

\newtheorem{theorem}{Theorem}

\newtheorem{lemma}[theorem]{Lemma}

\newtheorem{remark}[theorem]{Remark}

\newcommand{\eqa}{\begin{eqnarray}}
\newcommand{\eeqa}{\end{eqnarray}}
\newcommand{\beq}{\begin{equation}}
\newcommand{\eeq}{\end{equation}}
\newcommand{\nn}{\nonumber}

\newcommand{\pal}{\partial}

\newcommand{\pf}{\noindent{\it Proof \ }}

\newcommand{\epf}{$\quad$\hfill
\raisebox{0.11truecm}{\fbox{}}\par\vskip0.4truecm}

\title{On universality of critical behaviour in the focusing nonlinear Schr\"odinger equation, elliptic umbilic catastrophe and the {\it tritronqu\'ee} solution to the Painlev\'e-I equation.}
\author{{B.Dubrovin$^*$, T.Grava$^*$, C.Klein$^{**}$}\\
{\small $^*$ SISSA, Trieste; $^{**}$ Max-Planck Institute, Leipzig }}
\begin{document}
\maketitle
\begin{abstract} We argue that the critical behaviour near the point 
    of ``gradient catastrophe" of the solution to the Cauchy problem 
    for the focusing nonlinear Schr\"odinger equation $ 
    i\epsilon\,\Psi_t +\frac{\epsilon^2}2\Psi_{xx}+ |\Psi|^2 \Psi 
    =0$, $\epsilon\ll 1$,  with analytic initial data of the form $\Psi(x,0;\epsilon) =A(x) \, e^{\frac{i}{\epsilon} S(x)}$ is approximately described by a particular solution to the Painlev\'e-I equation.
\end{abstract}

\setcounter{equation}{0}
\setcounter{theorem}{0}
\section{Introduction}\par

The focusing nonlinear Schr\"odinger (NLS) equation for the complex valued function $\Psi=\Psi(x,t)$ 
\beq\label{i-nls0}
i\, \Psi_t +\frac{1}2\Psi_{xx}+ |\Psi|^2 \Psi =0
\eeq
has numerous physical applications in the description of nonlinear waves (see, e.g., the books \cite{whi, new, novik}). It can be considered as an infinite dimensional analogue of a completely integrable Hamiltonian system \cite{za}, where the Hamiltonian and the Poisson bracket is given by
\eqa\label{ham-nls}
&&
\Psi_t +\{ \Psi(x), {\cal H}\} =0
\nn\\
&&
\nn\\
&&
\{ \Psi(x), \Psi^*(y)\} =i\, \delta(x-y)
\\
&&
\nn\\
&&
{\cal H}=\frac12\,\int \left(  | \Psi_x |^2  - | \Psi\ |^4\right)\, dx
\nn
\eeqa
(here $\Psi^*$ stands for the complex conjugate function).
Properties of various classes of solutions to this equation have been 
extensively studied both analytically and numerically \cite{bronski, 
venak, ca, ce, fo, gr,   kam1, kam2, numart1d, miller, mi, to1, to2}. 
One of the striking features that distinguishes this equation from, 
say, the defocusing case
$$
i\, \Psi_t +\frac{1}2\Psi_{xx}- |\Psi|^2 \Psi =0
$$
is the phenomenon of {\it modulation instability} \cite{ag, ch, fo}. Namely, slow modulations of the plane wave solutions
$$
\Psi=A \, e^{i(kx-\omega t)},  \quad \omega=\frac12\, k^2-A^2
$$
develop fast oscillations in finite time.

The appropriate mathematical framework for studying these phenomena is the theory of the initial value problem
\beq\label{cauchy1}
\Psi(x,0;\epsilon) =A(x) \, e^{\frac{i}{\epsilon} S(x)}
\eeq
for the $\epsilon$-dependent NLS
\beq\label{nls0}
i\,\epsilon\, \Psi_t +\frac{\epsilon^2}2\Psi_{xx}+ |\Psi|^2 \Psi =0.
\eeq
Here $\epsilon>0$ is a small parameter, $A(x)$ and $S(x)$ are real-valued smooth functions.
Introducing the slow variables
\beq\label{slow}
u=|\Psi|^2, \quad v =\frac{\epsilon}{2i} \left( \frac{\Psi_x}{\Psi} - \frac{\Psi_x^*}{\Psi^*}\right)
\eeq
the equation can be recast into the following system:
\eqa\label{nls1}
&&
 u_t +(u\, v)_x=0
 \nn\\
 &&
 \\
 &&
v_t + v\, v_x - u_x +\frac{\epsilon^2}4 \left( \frac12 \frac{u_x^2}{u^2} - \frac{u_{xx}}{u}\right)_x=0.
\nn
\eeqa
The initial data for the system \eqref{nls1} coming from
\eqref{cauchy1}
do not depend on $\epsilon$:
\beq\label{cauchy2}
u(x,0)=A^2(x), \quad v(x,0)=S'(x).
\eeq
The simplest explanation of the modulation instability then comes from considering the so-called {\it dispersionless limit} $\epsilon \to 0$.
In this limit  one obtains the following first order quasilinear system
\beq\label{nls-d}
\left.\begin{array}{lr} u_t & + v\, u_x + u\, v_x=0 \\  & \\
v_t & -~~u_x + v\, v_x=0\end{array}\right\}.
\eeq
This is a system of {\it elliptic type} because of the condition $u>0$. Indeed, the eigenvalues of the coefficient matrix
$$
\left( \begin{array}{cc} v & u \\ -1 & v\end{array}
\right)
$$
are complex conjugate, $\lambda= v \pm i\sqrt{u}$. So, the Cauchy problem for the system \eqref{nls-d} is ill-posed in the Hadamard sense (cf. \cite{me, ca}). Even for analytic initial data the life span of a typical solution is finite, $t<t_0$. The $x$- and $t$-derivatives explode at some point $x=x_0$ when the time approaches $t_0$. This phenomenon is similar to the gradient catastrophe of solutions to nonlinear hyperbolic PDEs \cite{al}.

For the full system \eqref{nls1} the Cauchy problem is well-posed for a suitable class of $\epsilon$-independent initial data (see details in \cite{gwp1, gwp2}).  However, the well-posedness is not uniform in $\epsilon$. In practical terms that means that the solution to \eqref{nls1} behaves in a very irregular way in some regions of the $(x,t)$-plane when $\epsilon\to 0$. Such an irregular behaviour begins near the points $(x=x_0, t=t_0)$ of the ``gradient catastrophe" of the solution to the dispersionless limit \eqref{nls-d}. The solutions to \eqref{nls-d} and \eqref{nls1} are essentially indistinguishable for $t<t_0$; the situation changes dramatically near $x_0$ when approaching the critical point. Namely, when approaching $t=t_0$ the peak near a local maximum\footnote{Regarding initial data with local minima 
 we did not observe cusps related to minima in numerical
simulations. We believe they do not exist  because of   the focusing effect in the NLS that pushes maxima to cusps but seems to  smoothen
minima.}  of $u$ becomes more and more narrow due to self-focusing; 
the solution develops a zone of rapid oscillations for $t> t_0$. They 
have been studied both analytically and numerically in \cite{ce, fo, 
gr, jin, kam1, kam2, mi, to1, to2}. However, no results are available 
so far about the behaviour of the solutions to the focusing NLS {\it at the critical point} $(x_0, t_0)$.

The main subject of this work is the study of the behaviour of solutions to the Cauchy problem \eqref{nls1}, \eqref{cauchy2}
near the point of gradient catastrophe of the dispersionless system \eqref{nls-d}. In order to deal with the Cauchy problem for \eqref{nls-d} we will assume analyticity\footnote{We believe that the main conclusions of this paper must hold true also for non analytic initial data; the numerical experiments of \cite{ce} do not show much difference in the properties of solution between analytic and non analytic cases. However, the precise formulation of our Main Conjecture has to be refined in the non analytic case.} of the initial data $u(x,0)$, $v(x,0)$. Then the Cauchy problem for \eqref{nls-d} can be solved for $t<t_0$ via a suitable version of the hodograph transform (see Section \ref{sec1} below). An important feature of the gradient catastrophe for this system is that it happens at an {\it isolated point} of the $(x,t)$-plane, unlikely the case of KdV or defocusing NLS where the singularity of the hodograph solution takes place on a curve. We identify this singularity for a generic solution to \eqref{nls-d} as the {\it elliptic umbilic singularity} (see Section \ref{sec3} below) in the terminology of R.Thom \cite{th}. This codimension 2 singularity is one of the real forms labeled by the root system of the $D_4$ type in the terminology of V.Arnold {\sl et al.} \cite{ar}.

Our main goal is to find a replacement for the elliptic umbilic singularity when the dispersive terms are added, i.e., we want to describe the leading term of the asymptotic behaviour for $\epsilon\to 0$ of the solution to \eqref{nls1} near the critical point $(x_0, t_0)$ of a generic solution to \eqref{nls-d}. 

Thus, our study can be considered as a continuation of the programme 
initiated in \cite{du2} to study critical behaviour of Hamiltonian perturbations of nonlinear hyperbolic PDEs; the fundamental difference is that the non perturbed system \eqref{nls-d} is not hyperbolic! However, many ideas and methods of \cite{du2} (see also \cite{du1}) play an important role in our considerations.

The most important of these is the idea of {\it universality} of the 
critical behaviour. The general formulation of the universality 
suggested in \cite{du2} for the case of Hamiltonian perturbations of 
the scalar nonlinear transport equation says that the leading term of 
the multiscale asymptotics of the generic solution near the critical point does not depend on the choice of the solution, modulo Galilean transformations and rescalings. This leading term was identified in \cite{du2} via  a particular solution to the fourth order analogue of the Painlev\'e-I equation (the so-called $P_1^2$ equation). The existence of the needed solution to $P_1^2$ has been rigorously established in \cite{cl}. Moreover, it was argued in \cite{du2} that this behaviour is essentially independent on the choice of the Hamiltonian perturbation. Some of these universality conjectures have been partially confirmed by numerical analysis carried out in \cite{gr}.

The main message of this paper is the formulation of the Universality Conjecture for the critical behaviour of the solutions to the focusing NLS. Our considerations suggest the description of the leading term in the asymptotic expansion of the solution to \eqref{nls1} near the critical point via a particular solution to the classical Painlev\'e-I equation (P-I)
$$
\Omega_{\zeta\zeta} =6\, \Omega^2 -\zeta.
$$
The so-called {\it tritronqu\'ee} solution to P-I was discovered by 
P.Boutroux \cite{bo} as the unique solution having no poles in the sector $|\arg\zeta|<4\pi/5$ for sufficiently large $|\zeta|$. Remarkably, the very same solution\footnote{It is interesting that the same tritronqu\'ee solution (for real $\zeta$ only) appeares also in the study of certain critical phenomena in plasma \cite{sl}. In the theory of random matrices and orthogonal polynomials a different solution to P-I arises; see, e.g., \cite{duits}.}  arises in the critical behaviour of solutions to focusing NLS!

The paper is organized as follows. In Section \ref{sec1} we develop a 
version of the hodograph transform for integrating the dispersionless 
system \eqref{nls-d} before the catastrophe $t<t_0$. We also 
establish the shape of the singularity of the solution near the 
critical point; the latter is identified in Section \ref{sec3} with 
the elliptic umbilic catastrophe of Thom. In Section \ref{sec2} we 
develop a method of constructing formal perturbative solutions to the 
full system \eqref{nls1} before the critical time. In Section 
\ref{sec-pain} we collect the necessary information about the 
tritronqu\'ee solution of P-I and formulate the Main Conjecture of 
this paper. Such a formulation relies on a much stronger property of 
the tritronqu\'ee solution: namely, we need this solution to be {\it 
pole-free} in the whole sector $|\arg \zeta|<4\pi/5$. Numerical 
evidence for the absence of poles in this sector is given in Section \ref{sec5}. In Section \ref{sec6} we analyze numerically the agreement between the critical behaviour of solutions to focusing NLS and its conjectural description in terms of the tritronqu\'ee solution restricted on certain lines in the complex $\zeta$-plane. In the final Section \ref{sec7} we give some additional remark and outline the programme of future research.

\vskip 0.5truecm
\noindent{\bf Acknowledgments.} The authors thank K.McLaughlin for a very instructive discussion. One of the authors (B.D.) thanks R.Conte for bringing his attention to the tritronqu\'ees solutions of P-I. The results of this paper have been presented by one of the authors (T.G.) at the Conference ``The Theory of Highly Oscillatory Problems'', Newton Institute, Cambridge, March 26 - 30, 2007. T.G. thanks A.Fokas and S.Venakides for the stimulating discussion after the talk.
The present work is
partially supported by the European Science Foundation Programme ``Methods of
Integrable Systems, Geometry, Applied Mathematics" (MISGAM), Marie Curie RTN ``European Network in Geometry, Mathematical Physics and Applications"  (ENIGMA). The work of B.D. and T.G. is also partially supported 
by Italian Ministry of Universities and Researches (MUR) research grant PRIN 2006
``Geometric methods in the theory of nonlinear waves and their applications". 

\setcounter{equation}{0}
\setcounter{theorem}{0}
\section{Dispersionless NLS, its solutions and critical behaviour}\label{sec1}\par

The equations \eqref{nls1} are a Hamiltonian system
$$
\left. \begin{array}{r} u_t +\{ u(x), H\}=0 \\ \\ v_t +\{ v(x), H\}=0\end{array}\right\}
$$ 
with respect to the Poisson bracket originated in \eqref{ham-nls}
\beq\label{bra}
\{ u(x), v(y)\} =\delta'(x-y),
\eeq
other brackets vanish, with the Hamiltonian
\beq\label{ham}
H=\int \left[ \frac12 \left(u\, v^2 - u^2 \right) +\frac{\epsilon^2}{8u} u_x^2\right]\, dx.
\eeq 

Let us first describe the general analytic solution to the dispersionless system \eqref{nls-d}.

\begin{lemma} Let $u_I(x)$, $v_I(x)$ be two real valued analytic functions of the real variable $x$ satisfying 
$$
\left( u_{I,x}\right)^2 + \left( v_{I,x}\right)^2\neq 0.
$$
Then the solution $u=u(x,t)$, $v=v(x,t)$ to the Cauchy problem
\beq\label{cauchy}
u(x,0)=u_I(x), \quad v(x,0)=v_I(x)
\eeq
for the system (\ref{nls-d}) for sufficiently small $t$ can be determined
from the following system
\beq\label{hodo}
\left.\begin{array}{l} x= v\, t + f_u\\ \\ 0 = u\, t + f_v\end{array}\right\}
\eeq
where $f=f(u,v)$ is an analytic solution to the following linear PDE:
\beq\label{eq-f}
f_{vv}+u\, f_{uu}=0.
\eeq
Conversely, given any solution to (\ref{eq-f}) satisfying $u^2 f_{uu}^2 + f_{uv}^2 \neq 0$ at some point $(u=\tilde{u}, v=\tilde{v})$ such that $f_v(\tilde{u},\tilde{v})=0$, the system (\ref{hodo}) determines a solution to (\ref{nls-d}) defined locally near the point $x=\tilde{x}:= f_u(\tilde{u}, \tilde{v})$ for sufficiently small $t$.
\end{lemma}

\begin{remark} The solutions to the linear PDE (\ref{eq-f}) correspond to the first integrals of dispersionless NLS:
\beq\label{first-0}
F= \int f(u,v)\, dx, \quad \frac{d}{dt} F=0.
\eeq
Taking them as the Hamiltonians
\beq\label{flows}
\left.\begin{array}{r} u_s +\{u(x), F\} \equiv u_t + \left( f_v\right)_x=0 \\ \\
v_s +\{v(x), F\} \equiv v_t + \left( f_u\right)_x=0\end{array}\right\}
\eeq
yields infinitesimal symmetries of the dispersionless NLS:
\beq\label{comm}
\left( u_t\right)_s = \left( u_s\right)_t, \quad \left( v_t\right)_s = \left( v_s\right)_t.
\eeq
\end{remark}

One of the first integrals will be extensively used in this paper. It corresponds to the Hamiltonian density
\beq\label{ham-toda0}
g=-\frac12 v^2 + u(\log u -1).
\eeq
The associated Hamiltonian flow reads
\beq\label{toda-d}
\left.\begin{array}{r} u_s +v_x=0\\ \\ v_s = \frac{u_x}{u}\end{array}\right\}
\eeq
Eliminating the dependent variable $v$ one arrives at the elliptic version of the long wave limit of Toda lattice:
$$
u_{ss} + \left( \log u\right)_{xx}=0.
$$
Due to commutativity (\ref{comm}) the systems (\ref{nls-d}) and (\ref{toda-d}) admit a simultaneous solution $u=u(x,t,s)$, $v=v(x,t,s)$. Any such solution  can be locally determined from a system similar to (\ref{hodo})
\beq\label{hodo1}
\left.\begin{array}{l} x= v\, t + f_u\\ \\ s=u\, t + f_v\end{array}\right\}
\eeq
where $f=f(u,v)$, as above, solves the linear PDE (\ref{eq-f}).

The system (\ref{hodo1}) determines a solution $u=u(x,t,s)$, 
$v=v(x,t,s)$ provided applicability of the implicit function theorem. 
The conditions of the latter fail to hold at the {\it critical point}
$(x_0, s_0, t_0, u_0, v_0)$ such that
\beq\label{crit}
\begin{array}{l} x_0 =v_0 t_0 +f_u(u_0, v_0)\\ \\
s_0 = u_0 t_0 +f_v(u_0, v_0)\\ \\
f_{uu}(u_0, v_0)=f_{vv}(u_0, v_0)=0, \quad f_{uv}(u_0, v_0)=-t_0\end{array}
\eeq
In sequel we adopt the following system of notations: the values of the function $f$ and of its derivatives at the critical point will be denoted by $f^0$ {\sl etc}.   E.g., the last line of the conditions (\ref{crit}) will read
$$
f_{uu}^0 =f_{vv}^0=0, \quad f_{uv}^0 =-t_0.
$$

{\bf Definition 1}. We say that the critical point is {\it generic} if  at this point:
$$
f_{uuv}^0\neq 0.
$$
Let us the introduce real parameters $r$, $\psi$ determined by the third derivatives of the function $f$ evaluated at the critical point,
\beq\label{rpsi0}
\frac1{r} (\cos\psi -i\sin\psi) =  f_{uuv}^0 +i\sqrt{u_0} f_{uuu}^0.
\eeq
Due to the genericity assumption
$$
\psi\neq \frac{\pi}2 + \pi k.
$$

In order to describe the local behaviour of a solution to the dispersionless NLS/Toda equations we define a function $R(X; S, \psi)$ of real variables $X$, $S$ depending on the real parameter  $\psi$ satisfying
\beq\label{restrict}
(S+\cos\psi)^2 +(X+\sin\psi)^2 \neq 0
\eeq
by the following formula
\eqa\label{koren}
&&
R(X,S,\psi)=\sign[\cos{\psi}]
\\
&&
\times\sqrt{1+X \sin{\psi} +S\cos{\psi}+\sqrt{1+2 (X\sin{\psi} +S\cos{\psi}) +X^2 +S^2}}
\nn
\eeqa
Put
\eqa\label{koren1}
&&
P_0(X,S,\psi) = \frac1{\sqrt{2}} \left[ R(X,S,\psi)\cos{\psi} -
\frac{(X\cos{\psi}-S\sin{\psi})\sin{\psi}}{R(X,S,\psi)}\right] -\cos{\psi}
\nn\\
&&
\\
&&
Q_0(X,S,\psi)=\frac1{\sqrt{2}} \left[ \frac{(X\cos{\psi}-S\sin{\psi})\cos{\psi}}{R(X,S,\psi)} +R(X,S,\psi)\sin{\psi} \right] -\sin{\psi}.
\nn
\eeqa
Observe that $P_0(X, S,\psi)$ and $Q_0(X, S, \psi)$ are smooth functions of the real variable $X$ provided validity of the inequality (\ref{restrict}).

\begin{lemma} Given an analytic solution $u(x,s,t)$, $v(x,s,t)$  to the dispersionless NLS/Toda equations with a generic critical point $(x_0, s_0, t_0, u_0, v_0)$, and arbitrary real numbers $X$, $S$ satisfying (\ref{restrict}), $T<0$, then there exist the following limits
\eqa\label{limit-0}
&&
\lim_{\lambda\to +0} \lambda^{-1/2}\left[u(x_0 +\lambda^{1/2} v_0 T +\frac{\lambda}{2  \sqrt{u_0}} r\,X\, T^2, s_0 +\lambda^{1/2} u_0 T +\frac{\lambda}2 r S\, T^2, \lambda^{1/2} T)-u_0\right] 
\nn\\
&&
\nn\\
&&
=
r\, T\, P_0\left(X, S,\psi\right)
\nn\\
&&
\\
&&
\lim_{\lambda\to +0} \lambda^{-1/2}\left[v(x_0 +\lambda^{1/2} v_0 T +\frac{\lambda}{2  \sqrt{u_0}} r\, X\, T^2, s_0 +\lambda^{1/2} u_0 T +\frac{\lambda}2 r S\, T^2, \lambda^{1/2} T)-v_0\right] 
\nn\\
&&
\nn\\
&&
=
\frac{r}{\sqrt{u_0}}\, T\, Q_0\left(X, S,\psi\right)
\nn
\eeqa
where the parameters $r$, $\psi$ are defined by (\ref{rpsi0}).
\end{lemma}

\pf From the linear PDE (\ref{eq-f}) it follows that
$$
f_{uvv}= - u f_{uuu}-f_{uu}, \quad f_{vvv}=-u f_{uuv}.
$$
Using these formulae we expand the implicit function equations (\ref{hodo1}) near the critical point in the form
\eqa\label{umbilic}
&&
\bar x - v_0 \bar t = \bar v \, \bar t + \frac12 \left[ f_{uuu}^0 (\bar u^2 - u_0 \bar v^2) + 2 f_{uuv}^0 \bar u\, \bar v\right] +O\left((|\bar u|^2 +|\bar v|^2)^{3/2}\right)
\nn\\
&&
\\
&&
\bar s - u_0 \bar t = \bar u \, \bar t +\frac12 \left[ f_{uuv}^0  (\bar u^2 - u_0 \bar v^2)  -2 u_0 f_{uuu}^0 \bar u\, \bar v\right]+O\left((|\bar u|^2 +|\bar v|^2)^{3/2}\right)
\nn
\eeqa
where we introduce the shifted variables
\eqa
&&
\bar x=x-x_0, \quad \bar s=s-s_0, \quad \bar t=t-t_0
\nn\\
&&
\bar u=u-u_0, \quad \bar v=v-v_0.
\nn
\eeqa
The rescaling
\beq\label{scala}
\begin{array}{lcr}
\bar x - v_0 \bar t  & \mapsto &\lambda (\bar x - v_0 \bar t )\\
 & \\
\bar s - u_0 \bar t & \mapsto & \lambda (\bar s - u_0 \bar t )\\
 & \\
\bar t &\mapsto & \lambda^{1/2} \bar t\\
 & \\
\bar u & \mapsto & \lambda^{1/2} \bar u\\
 & \\
\bar v & \mapsto & \lambda^{1/2} \bar v
\end{array}
\eeq
in the limit $\lambda\to 0$ yields the quadratic equation
\beq\label{quadro}
z= \bar t\,w +\frac12 a\,w^2, \quad \bar t\neq 0
\eeq
where the complex independent and dependent variables $z$ and $w$ read
\beq\label{new}
z=\bar s+i\sqrt{u_0} \bar x - (u_0 +i\sqrt{u_0} v_0)\bar t, \quad w = \bar u+i\sqrt{u_0} \bar v
\eeq
and the complex constant $a$ is defined by
\beq\label{aa}
a= f_{uuv}^0 +i\sqrt{u_0} f_{uuu}^0,
\eeq
therefore
$$
\frac1{a} = r\, e^{i\psi}.
$$
The substitution
$$
X=2\, \sqrt{u_0}\,\frac{\bar x - v_0 \bar t}{ r\, \bar t^2 }, \quad S=2\,\frac{\bar s - u_0 \bar t}{ r\, \bar t^2 },
$$
reduces the quadratic equation to
$$
\left( w + \bar t\, r\, e^{i\psi}\right)^2 =r^2 \bar t^2 e^{2i\psi} \left[ 1+e^{-i\psi} (S+i\, X)\right].
$$
For $\bar t<0$ we choose the following root
\beq\label{root0}
w=r\,\bar t e^{i\,\psi} \left[ \sqrt{1+e^{-i\,\psi}(S+i\,X)}-1\right]
\eeq
where the branch of the square root is obtained by the analytic continuation of the one taking positive values on the positive real axis. Equivalently, 
$$
w =  \bar t \, r\,e^{i\psi}\left[ \sign(\cos{\psi})\frac{1}{\sqrt{2}} \left( \sqrt{\Delta +1 + S \cos\psi + X \sin\psi}+i \frac{X\cos{\psi} -S\sin{\psi}}{\sqrt{\Delta +1 + S \cos\psi + X \sin\psi}}\right)-1\right]
$$
where
$$
\Delta=\sqrt{1+2 (S\cos{\psi}+X\sin{\psi}) +X^2 +S^2}.
$$
This gives the formulae (\ref{koren}).  \epf

The result of the lemma describes the local structure of generic solutions to the dispersionless NLS/Toda equations near the critical point. It
can also be  represented in the following form
\eqa\label{asy}
&&
u(x,s,t) \simeq u_0 + r\, T\, P_0(X,S,\psi)
\nn\\
&&
\\
&&
v(x,s,t) \simeq v_0 + \frac1{\sqrt{u_0}} \, r\, T \, Q_0(X,S,\psi)
\nn
\eeqa
where
\beq\label{XST}
X=2\, \sqrt{u_0}\,\frac{\bar x - v_0 \bar t}{ r\, \bar t^2 }, \quad S=2\,\frac{\bar s - u_0 \bar t}{ r\, \bar t^2 }, \quad T=\bar t.
\eeq

We want to emphasize that the approximation \eqref{asy} works only near the critical point. Indeed, for large $x\to \pm\infty$ the function $u(x,s,t)$ and $v(x,s,t)$ have the following behaviour
\eqa\label{as-u}
&&
u= -\sqrt{r\, | x|}\, u_0^{1/4} \sqrt{1\mp\sin{\psi}} +u_0-r\, \bar t \cos{\psi} +O\left( \frac1{\sqrt{|x|}}\right)
\\
&&
\nn\\
&&
\label{as-v}
v =\mp\frac{\sqrt{r\, |x|}}{{u_0}^{1/4}} \sign(\cos{\psi}) \sqrt{1\pm\sin{\psi}}+v_0 -\frac{r}{\sqrt{u_0}}\, \bar t\, \sin{\psi}+O\left( \frac1{\sqrt{|x|}}\right).
\eeqa
So, for sufficiently large $|x|$ the function $u(x,s,t)$ defined by \eqref{asy} becomes negative.

The function $u$ has a maximum at the point $X=S\, \tan{\psi}$, so locally
\beq\label{ineq}
u\leq u_0 +r\, T \cos{\psi} -\sqrt{r}\, |\cos{\psi}| \, \sqrt{2\frac{S}{\cos{\psi}} +r\, T^2}<u_0.
\eeq
At the critical point $(x_0, s_0, t_0, u_0, v_0)$ the function $u$ develops a cusp. Let us consider only the particular case $S=0$ in order to avoid complicated expressions. In this case the local behaviour of the function $u$ near the critical point is given by
\beq\label{local}
\lim_{\bar t\to - 0}u=\left\{ \begin{array}{lr}u_0 -\sqrt{r\, |\hat x|} \,\sqrt{1-\sin{\psi}}, &\quad \hat x >0 \\ \\ u_0 -\sqrt{r\, |\hat x|}\, \sqrt{1+\sin{\psi}}, & \quad \hat x<0\end{array}\right. 
\eeq
(here $\hat x=\sqrt{u}_{0}(\bar x - v_0\bar t)$). Thus the parameters $r$, $\psi$ describe the shape of the cusp at the critical point. 

\setcounter{equation}{0}
\setcounter{theorem}{0}
\section{First integrals and solutions of the NLS/Toda equations}\label{sec2}\par

Let us first show that any first integral (\ref{first-0}) of the dispersionless equations can be uniquely extended to a first integral of the full equations.

\begin{lemma} Given a solution $f=f(u,v)$ to the linear PDE (\ref{eq-f}), there exists a unique, up to a total derivative, formal power series in $\epsilon$
$$
h_f=f+\sum_{k\geq 1} \epsilon^{2k} h_f^{[k]}(u, v; u_x, v_x, \dots, u^{(2k)}, v^{(2k)})
$$
such that the integral
$$
H_f =\int h_f\, dx
$$
commutes with the Hamiltonian of NLS equation:
$$
\{ H, H_f\}=0
$$
at every order in $\epsilon$. Explicitly,
\eqa\label{hf}
&&
h_f = f -\frac{\epsilon^2}{12} \left[ \left( f_{uuu} +\frac3{2u} f_{uu}\right) u_x^2 + 2 f_{uuv} u_x v_x -u f_{uuu} v_x^2\right]
\nn\\
&&
\\
&&
+\epsilon^4 \left\{ \frac1{120} \left[\left( f_{uuuu} +\frac5{2u} f_{uuu}\right) u_{xx}^2 +2 f_{uuuv} u_{xx} v_{xx} -u f_{uuuu} v_{xx}^2\right]\right.
\nn\\
&&
\nn\\
&&
- \frac1{80} f_{uuuu} u_{xx} v_x^2 -\frac1{48 u} f_{uuuv} v_{xx} u_x^2
-\frac1{3456 u^3} \left( 30 f_{uuu} -9 u f_{uuuu} + 12 u^2 f_{5 u} + 4 u^3 f_{6 u}\right) u_x^4
\nn\\
&&
\nn\\
&&
-\frac1{432 u^2} \left( -3 f_{uuuv} + 6 u f_{uuuuv} +2 u^2 f_{5u\, v}\right) u_x^3 v_x +\frac1{288 u} \left( 9 f_{uuuu} + 9 u f_{5 u} + 2 u^2 f_{6u}\right) u_x^2 v_x^2
\nn\\
&&
\nn\\
&&\left.
+\frac1{2160} \left( 9 f_{uuuuv} +10 u f_{5u\, v}\right) u_x v_x^3
-\frac{u}{4320}\left( 18 f_{5u} + 5 u f_{6u}\right) v_x^4\right\}+O(\epsilon^6)
\nn
\eeqa
\end{lemma}

Here we use short notations
$$
f_{5u} := \frac{\pal^5 f}{\pal u^5},\quad f_{6u} := \frac{\pal^6 f}{\pal u^6}, \quad f_{5u\, v} := \frac{\pal^6 f}{\pal u^5 \pal v}.
$$

{\bf Example 1}. Taking $f=\frac12 (u\, v^2 -u^2)$ one obtains the Hamiltonian of the NLS equation
$$
h_f = \frac12 (u\, v^2 -u^2)+\frac{\epsilon^2}{8u} u_x^2.
$$
In this case the infinite series truncates. It is easy to see that the series in $\epsilon$ truncates if and only if $f(u,v)$ is a polynomial in $u$.
Polynomial in $u$ solutions to the linear PDE (\ref{eq-f}) correspond to the standard first integrals of the NLS hierarchy.
\medskip

{\bf Example 2}. Taking $g=-\frac12 v^2 + u (\log u-1)$  (cf. (\ref{ham-toda0})) one obtains the Hamiltonian of Toda equation
\eqa
&&
h_g= -\frac12 v^2 + u (\log u-1)-\frac{\epsilon^2}{24 u^2} \left( u_x^2 + 2 u \, v_x^2\right) 
\nn\\
&&
\nn\\
&&
-\epsilon^4 \left( \frac{u_{xx}^2}{240 u^3} +\frac{v_{xx}^2}{60 u^2}
+\frac{u_{xx} v_{x}^2}{40 u^3} - \frac{u_x^4}{144 u^5} -\frac{u_x^2 v_x^2}{24 u^4} +\frac{v_x^4}{360 u^3}\right) +O(\epsilon^6)
\eeqa
written in terms of the function $\phi=\log u$ in the form
$$
\epsilon^2 \phi_{xx} +e^{\phi(s+\epsilon)} -2 e^{\phi(s)} +e^{\phi(s-\epsilon)}=0.
$$

\begin{lemma} Any solution to the NLS/Toda equations in the class of formal power series in $\epsilon$ can be obtained from the equations
\eqa\label{string}
&&
x=v\, t + \frac{\delta H_f}{\delta u(x)}
\nn\\
&&
\\
&&
s=u\, t + \frac{\delta H_f}{\delta v(x)}
\nn
\eeqa
where $f=f(u,v;\epsilon)$ is an arbitrary admissible solution to the linear PDE (\ref{eq-f}) in the class of formal power series in $\epsilon$,
$$
H_f=\int h_f\, dx.
$$
\end{lemma}

Now, we can apply to the system (\ref{string}) the rescaling (\ref{scala}) accompanied by the transformation
\beq\label{scala1}
\epsilon\mapsto \lambda^{5/4}\epsilon.
\eeq
At the limit $\lambda\to 0$ we arrive at the following system of equations
\eqa\label{3.5}
&&
\bar s - u_0 \bar t= \bar u\,\bar t +f_{uuv}^0 \left[ \frac12 (\bar u^2 -u_0 \bar v^2) +\frac{\epsilon^2}6 \bar u_{xx}\right] -u_0 f_{uuu}^0 \left[ \bar u \,\bar v +\frac{\epsilon^2}6 \bar v_{xx}\right]
\nn\\
&&
\\
&&
\bar x - v_0 \bar t = \bar v \,\bar t +f_{uuu}^0 \left[ \frac12 (\bar u^2 - u_0 \bar v^2) +\frac{\epsilon^2}6 \bar u_{xx}\right] +f_{uuv}^0 \left[ \bar u\,\bar v+\frac{\epsilon^2}6\bar v_{xx}\right].
\nn
\eeqa
Using the complex variables $z$, $w$ defined in (\ref{new}) we can rewrite the system in the following form:
\beq\label{pain1}
z\, r e^{i\psi} =w\, \bar t \,r e^{i\psi} +\frac12 w^2 +\frac{\epsilon^2}6 w_{xx}.
\eeq
The last observation is that the Toda equations generated by the Hamiltonian $H_g=\int h_g\, dx$ (see Example 2 above) after the scaling limit (\ref{scala}), (\ref{scala1}) yield the Cauchy - Riemann equations for the function $w=w(z)$,
$$
\pal w/\pal \bar z=0.
$$
Therefore the system \eqref{3.5} can be recast into the form equivalent to the Painlev\'e-I (P-I) equation (see \eqref{p-1} below)
\beq\label{pain2}
z\, r e^{i\psi} =w\, \bar t \,r e^{i\psi} +\frac12 w^2 -\frac{\epsilon^2}6 u_0 w_{zz}.
\eeq
Choosing
$$
\lambda=\epsilon^{4/5}
$$
we eliminate $\epsilon$ from the equation. 

In the Section \ref{sec-pain} below we will write explicitly the reduction of \eqref{pain2} to the Painlev\'e-I equation and give a conjectural characterization of the particular solution  of the latter.

\setcounter{equation}{0}
\setcounter{theorem}{0}
\section{Critical behaviour and elliptic umbilic catastrophe}\label{sec3}\par

Separating again the real and complex parts of  (\ref{pain1})
one obtains a system of  ODEs
\eqa\label{eqs}
&&
\frac{\epsilon^2}6 U_{XX} +\frac12 (U^2 -V^2) +r\,\bar t\, \left( U\cos{\psi} -V\sin{\psi}\right) -r\, \left( S\cos{\psi}-X\sin{\psi}\right)=0
\nn\\
&&
\\
&&
\frac{\epsilon^2}6 V_{XX} +UV +r\,\bar t\, \left( U\sin{\psi} +V\cos{\psi}\right) -r\, \left( S\sin{\psi} +X\cos{\psi}\right)=0
\nn
\eeqa
that can be identified with the Euler - Lagrange equations
$$
\delta S=0, \quad S=\int {\mathcal L}(U, V, U_X, V_X)\, dX
$$
with the Lagrangian
\eqa\label{lagr}
&&
{\mathcal L}=\frac{\epsilon^2}{12} \left( V_X^2 -U_X^2\right) +\frac16 \left( U^3 - 3 U\, V^2\right) +\frac12\,r\, \bar t\,  \left[ (U^2 - V^2) \cos\psi -2 U\, V\sin\psi\right] 
\nn\\
&&
\\
&&
+
r\, \left( X\sin\psi -S \cos\psi\right) U +r\, \left(S\sin\psi +X\cos\psi\right)V.
\nn
\eeqa
In the ``dispersionless limit" $\epsilon\to 0$ the Euler - Lagrange equations reduce to the search of stationary points of a function
(let us also set $\bar t=0$)
\beq\label{umb}
F=\frac16 \left( U^3 - 3 U\, V^2\right) + a_+ U + a_- V
\eeq
where we redenote
$$
a_+ = r\, \left( X\sin\psi -S \cos\psi\right), \quad a_-= r\, \left(S\sin\psi +X\cos\psi\right).
$$
At $a_+=a_-=0$ the function $F$ has an isolated singularity at the origin $U=V=0$ of the type $D_{4,\, -}$ also called {\it elliptic umbilic singularity}, according to R.Thom \cite{th}. This singularity appears in various physical problems; we mention here the caustics in the collisionless dark matter \cite{si} to give just an example.
The parameters $a_+$ and $a_-$ define two particular directions on the base of the miniversal unfolding of the elliptic umbilic; the full unfolding depending on 4 parameters reads
\beq\label{fold}
\hat F= \frac16 \left( U^3 - 3 U\, V^2\right) +\frac12 b(U^2 + V^2)+ a_+ U + a_- V +c.
\eeq
It would be interesting to study the properties of the modified Euler - Lagrange equations for the Lagrangian
$$
\hat{\mathcal L} = {\mathcal L}+\frac12 b(U^2 + V^2).
$$
This deformation does not seem to arrive from considering solutions to the NLS hierarchy.

\setcounter{equation}{0}
\setcounter{theorem}{0}
\section{The tritronqu\'ee solution to the Painlev\'e-I equation and the Main Conjecture}\label{sec-pain}\par

In this section we will select a particular solution to the Painlev\'e-I (P-I) equation
\beq\label{p-1}
\Omega_{\zeta\zeta} =6 \Omega^2 -\zeta
\eeq
Recall \cite{in} that an arbitrary solution to this equation is a meromorphic function on the complex $\zeta$-plane. According to P. Boutroux \cite{bo} the poles of the solutions accumulate along the rays
\beq\label{ray}
\arg \zeta = \frac{2\pi n}{5}, \quad n=0,\, \pm 1, \, \pm 2.
\eeq
Boutroux proved that, for each ray there is a one-parameter family of particular solutions called {\it int\'egrales tronqu\'ees} whose lines of poles truncate for large $\zeta$. He proved that the
{\it int\'egrale tronqu\'ee}  has no poles for large $|\zeta|$ within two consecutive sectors of the angle $2\pi/5$ around the ray, and, moreover it has the asymptotic behaviour of the form
\beq\label{as-bout}
\Omega =-\left(\frac{\zeta}6\right)^{1/2} \left[ 1 +O\left( \zeta^{-\frac34(1-\varepsilon)}\right)\right]
\eeq
for a suitable choice of the branch of the square root (see below) and a sufficiently small $\varepsilon >0$.

Furthemore, if a solution truncates along any two of the rays \eqref{ray} then it truncates along three of them. These particular solutions to P-I are called {\it tritronqu\'ees}. They have no poles for large $|\zeta|$ in {\it four} consecutive sectors; their asymptotics for large $\zeta$ is given by \eqref{as-bout}. It suffices to know the tritronqu\'ee solution $\Omega_0(\zeta)$ for the sector
\beq\label{sec}
|\arg \zeta | < \frac{4\pi}5.
\eeq
In this case the branch of the square root in \eqref{as-bout} is obtained by the analytic continuation of the principal branch taking positive values on the positive half axis $\zeta >0$. Other four tritronqu\'ees solutions are obtained by applying the symmetry
\begin{equation}
    \Omega_n(\zeta) =e^{\frac{4\pi i n}5} \Omega_0\left( e^{\frac{2\pi i n}5} \zeta\right), \quad n=\pm 1, \quad \pm 2.
    \label{triphi}
\end{equation}
The properties of the tritronqu\'ees solutions in the finite part of the complex plane were studied in the important paper of N.Joshi and A.Kitaev \cite{jk}.

A. Kapaev \cite{ka} obtained a complete characterization of the tritronqu\'ees solutions in terms of the Riemann - Hilbert problem associated with P-I. We will briefly sketch here the main steps of his construction.

The equation \eqref{p-1} can be represented as the compatibility condition of the following system of linear differential equations for a two-component vector valued function $\Phi=\Phi(\lambda,\zeta)$
\eqa\label{u-lam}
&&
\Phi_\lambda=\left(\begin{array}{cc} \Omega_\zeta & 2\lambda^2 + 2 \Omega\, \lambda -\zeta + 2 \Omega^2 \\
 & \\
 2(\lambda-\Omega) & -\Omega_\zeta\end{array}\right)\, \Phi
 \\
 &&
 \nn\\
 &&
 \nn\\
 &&
 \Phi_\zeta=- \left(\begin{array}{cc} 0 & \lambda+2 \Omega\\
  & \\ 
  1 & 0\end{array}\right)\, \Phi.
  \label{u-zeta}
  \eeqa
The canonical matrix solutions $\Phi_k(\lambda,\zeta)$ to the system \eqref{u-lam} - \eqref{u-zeta} are uniquely determined by their asymptotic behaviour
\eqa\label{asym}
&&
\Phi_k(\lambda,\zeta) \sim 
\frac1{\sqrt{2}} \left(\begin{array}{cc}\lambda^{1/4} & \lambda^{1/4}\\ \lambda^{-1/4} & -\lambda^{-1/4}\end{array}\right) 
\left[ 1 -\frac1{\sqrt{\lambda}}\left(\begin{array}{cc}H & 0 \\ 0 & -H\end{array}\right)\right.
\nn\\
&&
\\
&&\left.
\quad\quad\quad\quad\quad +\frac1{2\lambda}\left(\begin{array}{cc} H^2 & \Omega \\ \Omega & H^2\end{array}\right) +O\left( \lambda^{-3/2}\right) \right]\, e^{\theta(\lambda,\zeta)\, \sigma_3}, \quad |\lambda|\to\infty, \quad \lambda\in\Sigma_k
\nn
\eeqa
in the sectors
\beq\label{sec-k}
\Sigma_k=\left\{ \lambda\in\mathbb{C}\, | \, \frac{2\pi}{5}\left( k-\frac32\right) < \arg \lambda < \frac{2\pi}{5}\left( k+\frac12\right)\right\}, \quad k\in\mathbb{Z}.
\eeq
Here
\beq\label{theta}
\theta(\lambda,\zeta) =\frac45 \lambda^{5/2} -\zeta\, \lambda^{1/2}, \quad \sigma_3 =\left(\begin{array}{cc}1 & 0\\ 0 & -1\end{array}\right), \quad H = \frac12\, \Omega_\zeta^2 -2\,\Omega^3 + \zeta\, \Omega,
\eeq
the branch cut on the complex $\lambda$-plane for the fractional powers of $\lambda$ is chosen along the negative real half-line.

The {\it Stokes matrices} $S_k$ are defined by
\beq\label{stokes1}
\Phi_{k+1}(\lambda,\zeta) =\Phi_k(\lambda,\zeta) \, S_k, \quad \lambda\in \Sigma_k \cap \Sigma_{k+1}.
\eeq
They have the triangular form
\beq\label{stokes2}
S_{2k-1} =\left( \begin{array}{cc} 1 & s_{2k-1}\\ 0 & 1\end{array} \right), \quad S_{2k} =\left(\begin{array}{cc} 1 & 0\\ s_{2k} & 1\end{array}\right)
\eeq
and satisfy the constraints
\beq\label{stokes3}
S_{k+5} =\sigma_1\, S_k\, \sigma_1, \quad \quad k\in\mathbb{Z}; \quad S_1 S_2 S_3 S_4 S_5 =i\, \sigma_1 
\eeq
where
$$
\sigma_1=\left(\begin{array}{cc} 0 & 1 \\ 1 & 0\end{array}\right).
$$
Due to \eqref{stokes3} two of the {\it Stokes multipliers} $s_k$ 
determine all others; they depend neither on $\lambda$ nor on $\zeta$ provided $\Omega(\zeta)$ satisfies \eqref{p-1}. 

In order to obtain a parametrization of solutions to the P-I equation 
\eqref{p-1} by Stokes multipliers of the linear differential 
equation  \eqref{u-lam} one has to reformulate the above definitions as a certain Riemann - Hilbert problem. The solution of the Riemann - Hilbert problem depends on $\zeta$ through the asymptotics \eqref{asym}. If the Riemann - Hilbert problem has a unique solution for the given $\zeta_0\in\mathbb{C}$ then the canonical matrices $\Phi_k(\lambda,\zeta)$ depend analytically on $\zeta$ for sufficiently small $|\zeta-\zeta_0|$; the coefficient $\Omega=\Omega(\zeta)$ will then satisfy \eqref{p-1}. The {\it poles} of the meromorphic function $\Omega(\zeta)$ correspond to the forbidden values of the parameter $\zeta$ for which the Riemann - Hilbert problem admits no solution.

We will now consider a particular solution to the P-I equation specified by the following Riemann - Hilbert problem. Denote four oriented rays $\gamma_0$, $\gamma_{\pm 1}$, $\rho$ in the complex $\lambda$-plane defined by
\eqa\label{gam}
&&
\gamma_k=\{ \lambda\in\mathbb{C} \, | \, \arg \lambda=\frac{2\pi k}5\}, \quad k=0, \, \pm 1
\nn\\
&&
\\
&&
\rho=\{ \lambda\in\mathbb{C} \,\, ~ | \, \arg \lambda=\pi\}
\nn
\eeqa
directed towards infinity. The rays divide the complex plane in four sectors. We are looking for a piecewise analytic function
$\Phi(\lambda, \zeta)$ on 
$$
\lambda\in \mathbb{C}\setminus (\gamma_{-1}\cup \gamma_0 \cup \gamma_1 \cup \rho)
$$
depending on the parameter $\zeta$ continuous up to the boundary with the asymptotic behaviour at $|\lambda|\to \infty$ of the form \eqref{asym} satisfying the following jump conditions on the rays
\eqa\label{rh1}
&&
\Phi_+(\lambda, \zeta)  = \Phi_-(\lambda, \zeta) S_k, \quad \lambda\in \gamma_k
\nn\\
&&
\\
&&
\Phi_+ (\lambda, \zeta) = \Phi_-(\lambda, \zeta) S_\rho, \quad \lambda\in \rho.
\nn
\eeqa
Here the plus/minus subscripts refer to the boundary values of $\Phi$ respectively on the left/right sides of the corresponding oriented ray, the jump matrices are given by
\beq\label{rh2}
S_0= \left(\begin{array}{cc} 1 & 0\\ i & 1\end{array}\right), \quad S_{\pm 1} =\left(\begin{array}{cc} 1 & i\\ 0 & 1\end{array}\right), \quad S_\rho = \left( \begin{array}{cc} 0 & -i\\ -i & 0\end{array}\right).
\eeq

The following result is due to A.Kapaev\footnote{Our solution $\Omega_0(\zeta)$ coincides with $y_3(x)\equiv y_{-2}(x)$, $x=-\zeta$,  of \cite{ka} (see eq. (2.73) of \cite{ka}; Kapaev uses a different normalization $y''=6 y^2 +x$ of the P-I equation).} .

\begin{theorem}The solution to the above Riemann - Hilbert problem exists and it is unique for
\beq\label{sector}
|\arg\lambda | < \frac{4\pi}5, \quad |\lambda| >R
\eeq
for a sufficiently large positive number $R$. The associated function 
\eqa\label{resh}
&&
\Omega_0(\zeta) := \frac{dH(\zeta)}{d\zeta},
\nn\\
&&
\\
&& 
H(\zeta):= \left[ \lim_{\lambda\to\infty} \lambda^{1/2} \left( \frac1{\sqrt{2}} \left( \begin{array}{cc} 1 & 1\\ 1 & -1\end{array}\right) \, \lambda^{-\frac14\, \sigma_3} \Phi(\lambda,\zeta) \, e^{-\theta(\lambda,\zeta) \, \sigma_3} -1\right)\right]_{11}
\nn
\eeqa
is analytic in the domain \eqref{sector}, it satisfies P-I and enjoys the asymptotic behaviour
\beq\label {koren2}
\Omega_0(\zeta)\sim -\sqrt{\frac{\zeta}6}, \quad |\zeta| \to \infty, \quad |\arg\lambda | < \frac{4\pi}5.
\eeq
Moreover, any solution of P-I having no poles in the sector \eqref{sector} for some large $R>0$ coincides with $\Omega_0(\zeta)$.
\end{theorem}

Joshi and Kitaev proved that the {\it tritronqu\'ee} solution has no poles on the positive real axis. They found a numerical estimate for the position of the first pole $\zeta_0$ of the tritronqu\'ee solution $\Omega_0(\zeta)$ on the negative real axis:
$$
\zeta_0\simeq -2.3841687
$$
(cf. also \cite{co}).
Besides this estimate very little is known about the location of poles of $\Omega_0(\zeta)$. Our numerical experiments (see below) suggest the following

{\bf Main Conjecture. Part 1}. {\it The tritronqu\'ee solution $\Omega_0(\zeta)$ has no poles in the sector
\beq\label{sector1}
|\arg\lambda | < \frac{4\pi}5.
\eeq
}

\medskip

We are now ready to describe the conjectural universal structure behind the critical behaviour of generic solutions to the focusing NLS. For simplicity of the formulation let us assume $\cos\psi >0$.

\medskip

{\bf Main Conjecture. Part 2}. {\it Any generic solution to the NLS/Toda equations near the critical point behaves as follows
\eqa\label{main}
&&
u(x,s,t;\epsilon)+i\sqrt{u_0} v(x,s,t;\epsilon)
\simeq u_0 +i\sqrt{u_0}v_0 -\bar t\, r e^{i\psi} +2 \,\epsilon^{2/5}
(3 r\sqrt{u_0} )^{2/5} e^{\frac{2i\psi}5} \Omega_0(\zeta) +O\left( \epsilon^{4/5}\right) 
\nn\\
&&
\\
&&
\zeta=\left( \frac{3r}{u_0^2}\right)^{1/5} e^{\frac{i\psi}5} 
\left[ \frac{\bar s - u_0 \bar t +i\sqrt{u_0} (\bar x - v_0\bar t)+\frac12 re^{i\psi} \bar t^2}{\epsilon^{4/5}}\right]
\nn 
\eeqa
where $\Omega_0(\zeta)$ is the tritronqu\'ee solution to the  Painlev\'e-I equation \eqref{p-1}.
} 

\medskip

The above considerations can actually be applied replacing the NLS time flow by any other flow of the NLS/Toda hierarchy. The local description of the critical behaviour remains unchanged.
\begin{remark}
    Note that the angle of the line $\zeta(\bar{x})$ in 
    (\ref{main}) for $\bar{t}$ fixed is equal to $\psi/5+\pi/2$, 
    $\psi\in[-\pi,\pi]$. Thus 
    the maximal value of $\mbox{arg}\zeta$ is equal to $7\pi/10<4\pi/5$. The lines in 
    (\ref{main}) consequently do not get close to the critical lines 
    of the \emph{tritronqu\'ee} solution. 
\end{remark}

\setcounter{equation}{0}
\setcounter{theorem}{0}
\setcounter{equation}{0}
\setcounter{theorem}{0}
\section{Numerical analysis of the {\it tritronqu\'ee} solution of P-I}\label{sec5}\par

In this section we will numerically construct the \emph{tritronqu\'ee} 
solution $\Omega_{0}$, i.e.\ the \emph{tritronqu\'ee} solution with 
asymptotic behavior \eqref{as-bout}. We will drop the index $0$ in 
the following. 
The solution will be first constructed on a straight line in the complex plane. In a 
second step we will then explore global properties of these solutions 
within the limitations imposed by a numerical approach\footnote{Cf.~\cite{fokas}
where a similar technique was applied to solve  numerically the Painlev\'e-II equation in the complex domain.}.

Let the straight line in the complex plane be given by $\zeta=ay+b$ 
with $a,b\in\mathbb{C}$ constant (we choose $a$ to have a non-negative imaginary 
part) and $y\in\mathbb{R}$. The asymptotic 
conditions are 
\begin{equation}
    \Omega \sim -\sqrt{\frac{\zeta}{6}}
    \label{PI2},
\end{equation}
for $y\to\pm\infty$.
 The root is defined to have its cut  along the negative real axis and to assume positive values on the positive real axis. 
This choice of the root implies the following symmetry for the 
solution:
\begin{equation}
    \Omega(\zeta^{*}) = \Omega^{*}(\zeta)
    \label{PIsym}.
\end{equation}
Thus $\Omega$ is real on the real axis, see \cite{jk}.

Numerically it is not convenient to impose boundary conditions at 
infinity. We thus assume that the wanted solution  can be expanded  in a Laurent 
series in $\sqrt{\zeta}$ around infinity. Such an asymptotic expansion is 
possible for the considered \emph{tritronqu\'ee} solution in the sector 
$|\arg \zeta|<4\pi/5$. The formal series can be written there (see 
\cite{jk}) in the form 
\begin{equation}
    \Omega_{f}=-\sqrt{\frac{\zeta}{6}}\sum_{k=0}^{\infty}\frac{a_{k}}{\zeta^{5k/2}}
    \label{formal},
\end{equation}
where $a_{0}=1$, and where the remaining coefficients follow from the 
recurrence relation for $k\geq0$
\begin{equation}
    a_{k+1}= 
    \frac{25k^{2}-1}{8\sqrt{6}}a_{k}-\frac{1}{2}\sum_{m=1}^{k}a_{m}a_{k+1-m}.
    \label{forrec}
\end{equation}
This formal series is divergent, the coefficients $a_{k}$ behave 
asymptotically as $((k-1)!)^{2}$, see \cite{jk} for a detailed 
discussion. 

It is known that divergent series can be used to get numerically 
acceptable approximations to the values of the sum by properly 
truncating the series. Generally  the best approximations for the sum 
result from truncating the series where the terms take the smallest 
absolute values  (see e.g.~\cite{gradrhy}). Since we work in Matlab with a 
precision of 16 digits and with values of $|\zeta|\geq 10$, we 
typically consider up to 10 terms in the series. In this case the 
terms corresponding to $a_{10}$ are of the order of machine precision 
($10^{-14}$ and below).

Thus we have constructed approximations to the numerical values of 
the \emph{tritronqu\'ee} solution for large values of $|\zeta|$. These can 
be used as in \cite{jk} to set up an initial value problem for the P-I 
equation and to solve this with a standard ODE solver. 
In fact the approach works well on the real axis starting from 
positive values until one reaches the first singularity on the 
negative real axis. It is straightforward to check the results of 
\cite{jk} with e.g.~\emph{ode45}, the Runge-Kutta solver in Matlab 
corresponding to the Maple solver used in \cite{jk}. If one solves 
P-I on a line that avoids the sector $|\arg\zeta|>4\pi/5$, one 
could integrate until one reaches once more large values of $|\zeta|$ 
for which the asymptotic conditions are known. This would provide a 
control on the numerical accuracy of this so-called shooting 
approach. 
Shooting methods are problematic if the second solution to the 
initial value problem has poles as is the case for P-I. In this case the 
numerical errors in the initial data (here due to the asymptotic 
conditions) and in the time integration will lead to a large 
contribution of the unwanted solution close to its poles which will 
make the numerical solution useless. It is obvious that P-I has such 
poles from the numerical results in \cite{jk} and the property 
(\ref{triphi}). In \cite{jk} the task was to locate poles in the 
\emph{tritronqu\'ee} solution, and in this case the shooting 
approach seems to be the only available. Here we are studying, 
however, the solution on a line in the complex plane where we know the 
asymptotic conditions for the affine parameter tending to $\pm\infty$.

\begin{figure}[!htb]
\centering
\includegraphics[width=5in]{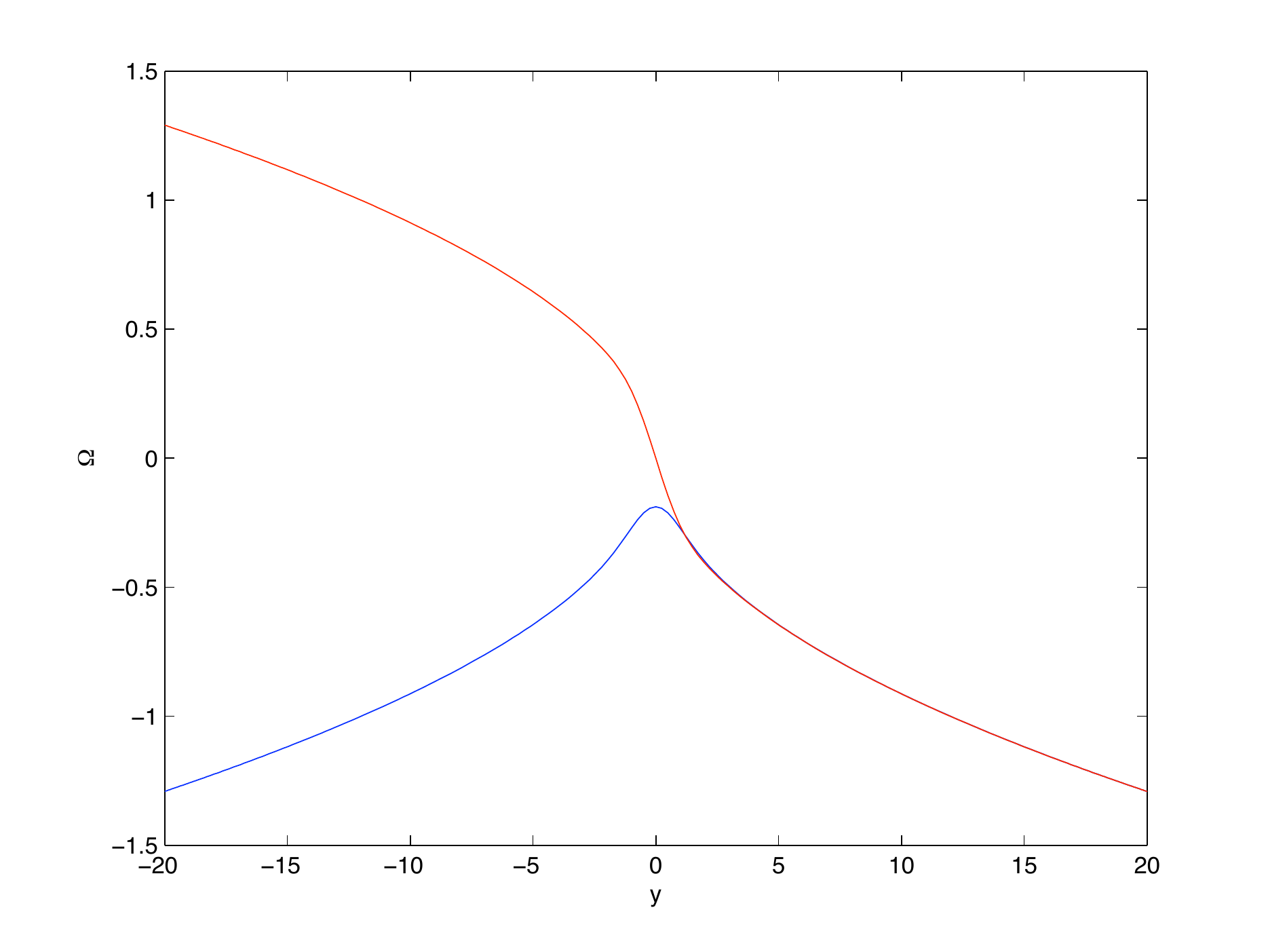}
\caption{Real (blue) and imaginary part (red) of the \emph{tritronqu\'ee}
solution to 
the Painlev\'e I equation for $\zeta=iy$.}
\label{figpain2}
\end{figure}

\begin{figure}[!htb]
\centering
\includegraphics[width=5in]{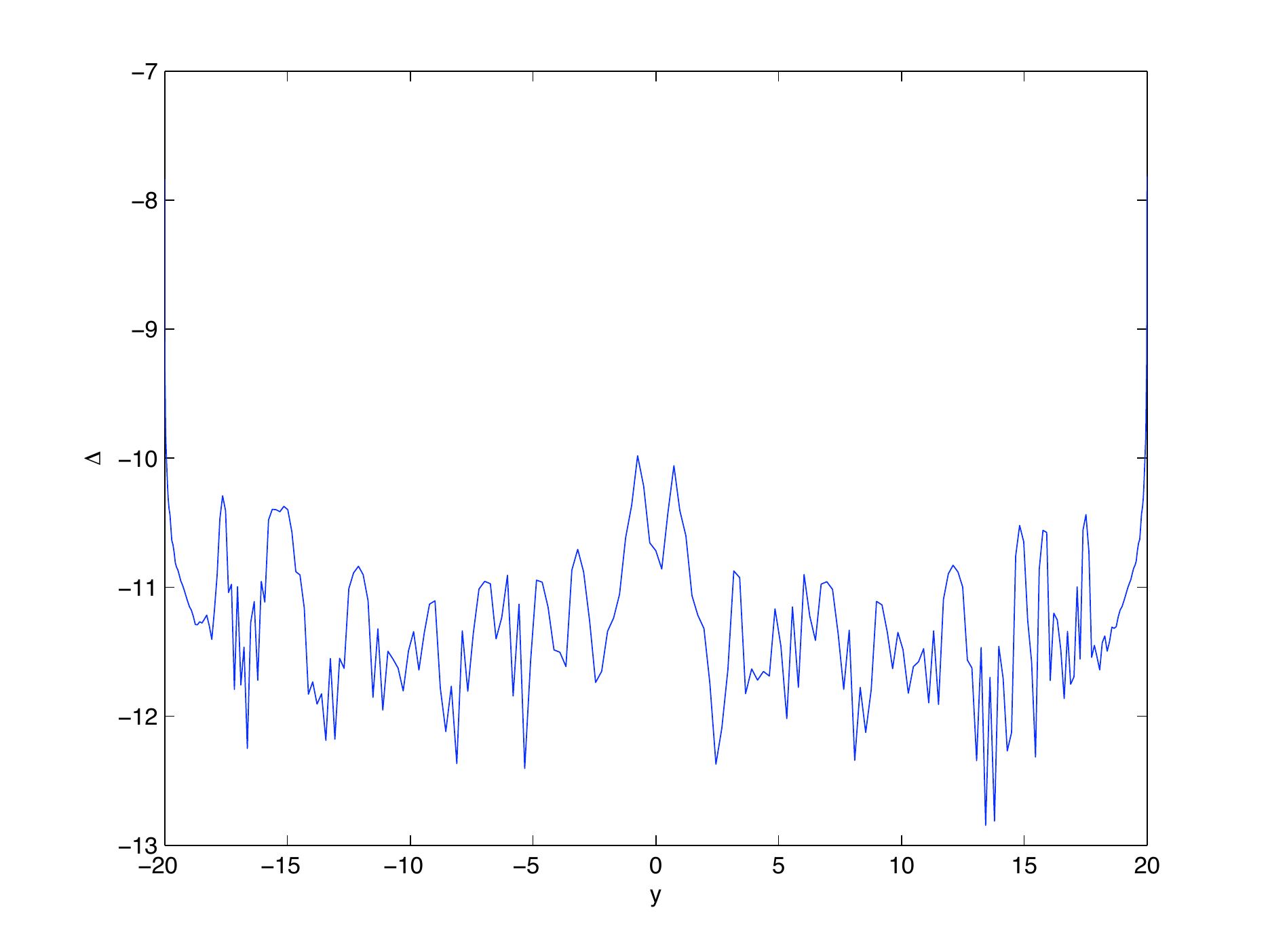}
\caption{Error in the solution of the Painlev\'e I equation.}
\label{figp2test}
\end{figure}

Thus we use as in \cite{fokas} the asymptotic conditions on lines 
avoiding the sector $|\arg\zeta|>4\pi/5$ to set up a boundary 
value problem for $y=\pm y_{0}$, $y_{0}\geq 10$. 
The solution in the interval $[-y_{0},y_{0}]$ is numerically obtained 
with a finite difference code based on a collocation method. The code 
\emph{bvp4c} distributed with Matlab, see \cite{bvp4c} for details, 
uses cubic polynomials in between the collocation points. The 
P-I equation is rewritten in the form of a first order 
system. With some 
initial guess (we use $\Omega=-\sqrt{\zeta/6}$ as the initial guess), the 
differential equation is solved iteratively by linearization.  
The collocation points (we use up to 10000) 
are dynamically adjusted during 
the iteration. The iteration is stopped when the equation is 
satisfied at the collocation points with a prescribed relative accuracy, 
typically $10^{-10}$. 
The solution for $a=i$ and $b=0$ is shown in Fig.~\ref{figpain2}.
The values of $\Omega$ in between the collocation 
points are obtained by interpolation via the cubic polynomials in terms of which the 
solution has been constructed. This interpolation leads to a loss in 
accuracy of roughly one order of magnitude with respect to the 
precision at the collocation points. To test this we determine 
the numerical solution via \emph{bvp4c} for P-I
on Chebychev collocation 
points and check the accuracy with which the equation is satisfied 
via Chebychev 
differentiation, see e.g.\ \cite{trefethen1}.  It is found 
that the numerical solution with a relative tolerance of $10^{-10}$ on 
the collocation points satisfies the ODE to roughly the same 
order except at the boundary points where it is of the order $10^{-8}$, 
see Fig.~\ref{figp2test} where we show the residual $\Delta$ by 
plugging the numerical solution into the differential equation for 
the above example. It is straightforward to achieve a prescribed
accuracy by requiring a certain value for the relative tolerance. 
Notice that we are not interested in a high precision solution of P-I 
here, but in a comparison of solutions to the NLS equation close to the 
point of gradient catastrophe of the semiclassical system with an 
asymptotic solution in terms of P-I transcendents.
For this purpose an accuracy of the solution of the order of 
$10^{-4}$ will be sufficient in all studied cases.  

\begin{figure}[!htb]
\centering
\includegraphics[width=5in]{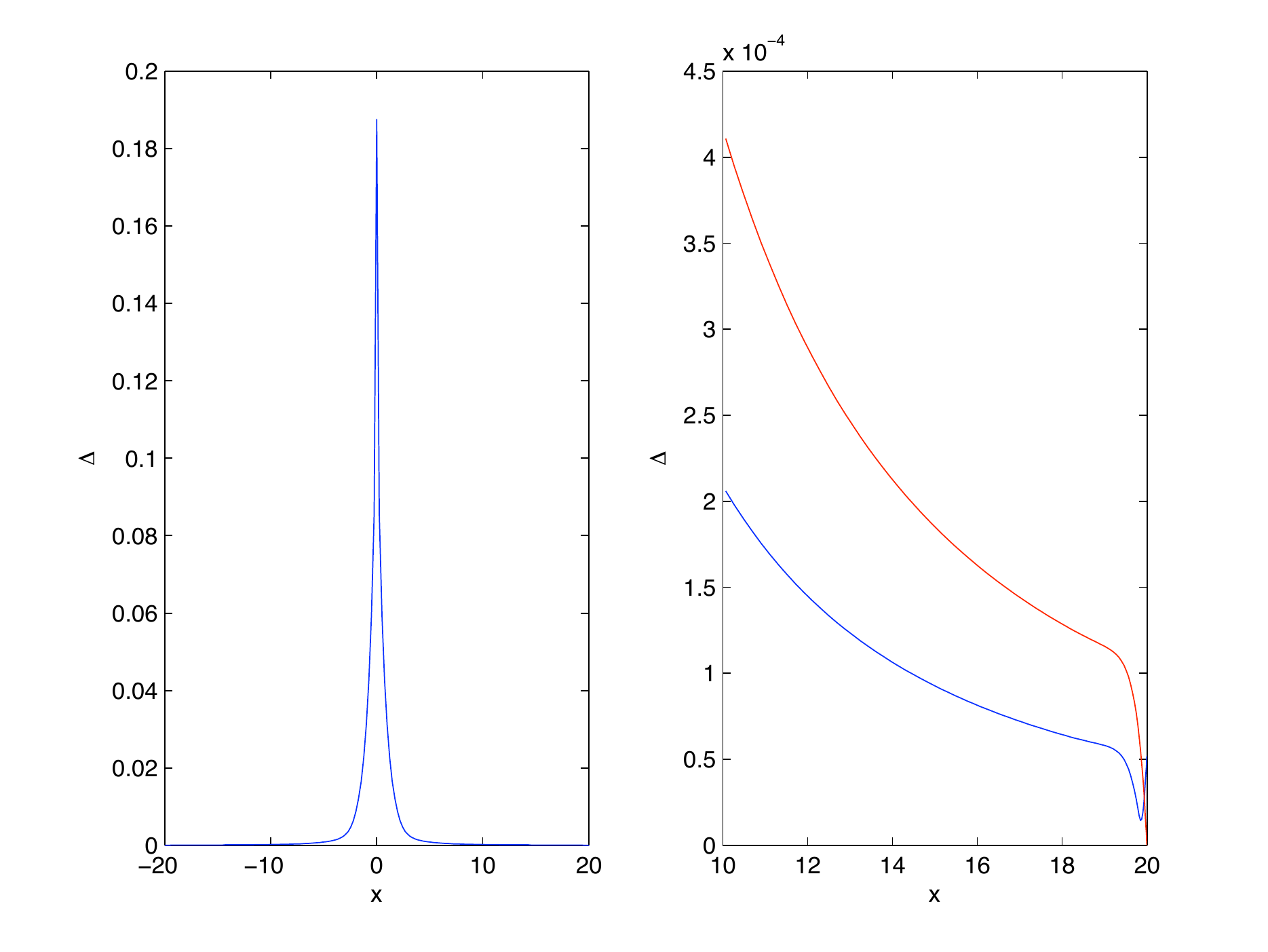}
\caption{The plot on the left side shows the absolute value of the 
difference between the \emph{tritronqu\'ee} solution and the asymptotic 
condition $-\sqrt{\zeta/6}$ for $a=i$ and $b=0$. The plot on the 
right side shows 
in blue the same difference for $x>10$ and in red the difference 
between the \emph{tritronqu\'ee} solution and the truncated asymptotic 
series.}
\label{figpain1_0asym}
\end{figure}

The quality of the used boundary conditions via the asymptotic 
behavior can be checked by computing the solution for 
different values of $y_{0}$. One finds that the difference between 
the asymptotic square root and the \emph{tritronqu\'ee} solution is only 
visible near the origin, see Fig.~\ref{figpain1_0asym}. For large $x$ 
it can be seen that the difference between the square root 
asymptotics and the \emph{tritronqu\'ee} solution reaches quickly values 
below the aimed at threshold of $10^{-4}$. It is interesting to note 
that this difference is actually smaller than the difference between 
the \emph{tritronqu\'ee} solution and the truncated formal asymptotic 
series except at the boundary, where the latter condition is implemented 
(see Fig.~\ref{figpain1_0asym}). 

\begin{figure}[!htb]
\centering
\includegraphics[width=5in]{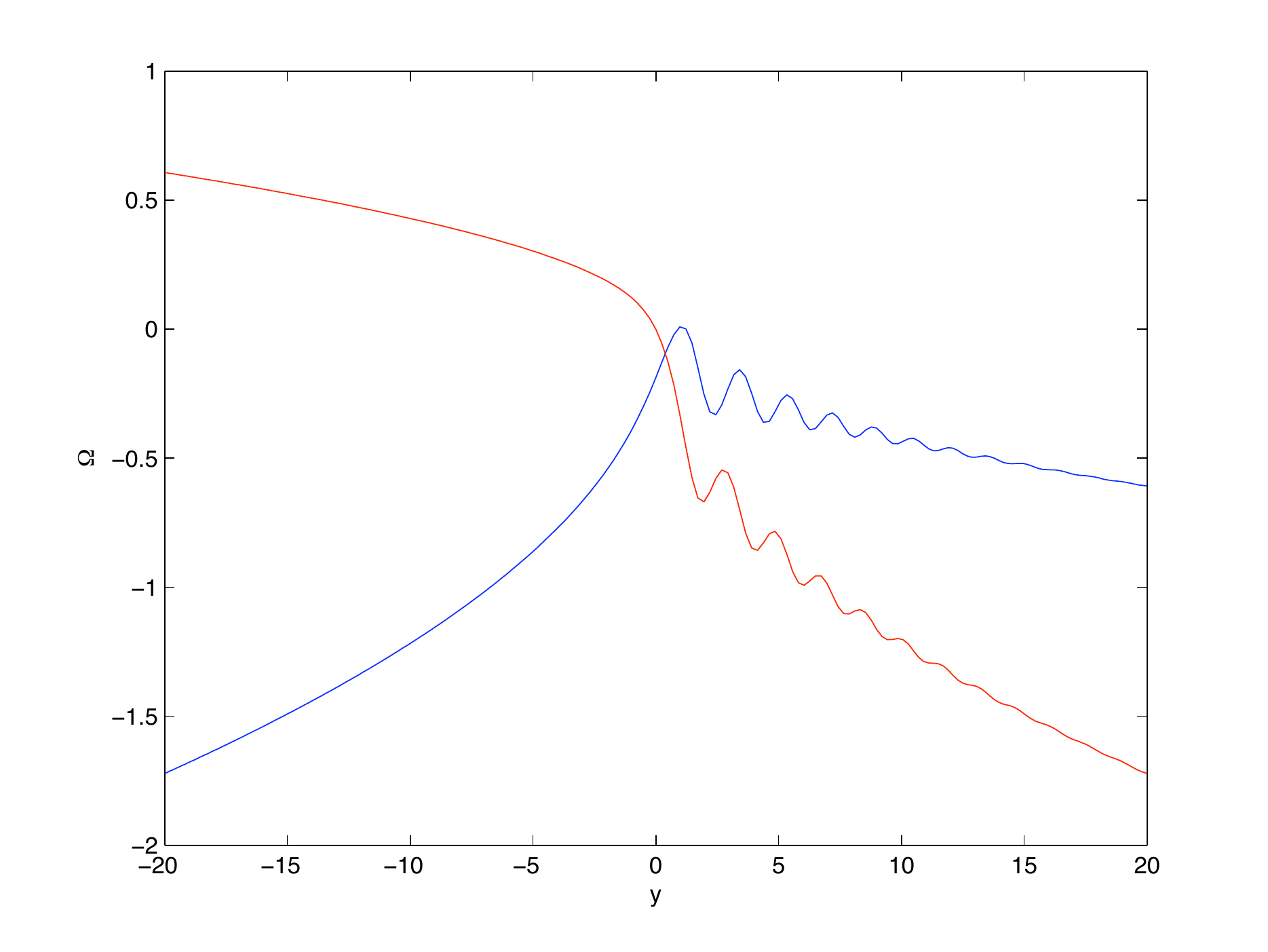}
\caption{Real (blue) and imaginary part (red) of the \emph{tritronqu\'ee} 
solution close to the critical line (for $a = 
\exp(i(4\pi /5-.05))$) with oscillations of slowly decreasing amplitude.}
\label{figOm_.05}
\end{figure}

The dominant behavior of the square root changes if one approaches 
the critical lines $a = \exp(4\pi i/5)$, $b=0$. As can be seen from 
Fig.~\ref{figOm_.05}, the solution shows oscillations on top of the 
square root. The closer one comes to the critical lines, the slower 
is the fall off of the amplitude of the oscillations. 
We conjecture that these oscillations will have on the critical lines 
only a slow algebraic fall off towards infinity. 

The above approach thus allows the computation of the \emph{tritronqu\'ee} 
solution for a line avoiding the sector $|\arg\zeta|>4\pi/5$ 
with high accuracy. The picture one obtains by computing $\Omega$ 
along several such lines is that there are indeed no singularities in 
the sector $|\arg\zeta|<4\pi/5$, and that the square root 
behavior is followed for large $|\zeta|$. To obtain a more complete 
picture, we compute the \emph{tritronqu\'ee} solution for 
$|\arg\zeta|<4\pi/5-0.05$ and $|\zeta|<R$ ( we choose $R=20$). 
The boundary data 
for $|\zeta|=R$ follow as before from the truncated asymptotic 
series, the data for $\arg\, \zeta=\pm 4\pi/5-0.05$ are obtained 
by computing the \emph{tritronqu\'ee} solution on the respective lines as 
above.

\begin{figure}[!htb]
\centering
\includegraphics[width=5in]{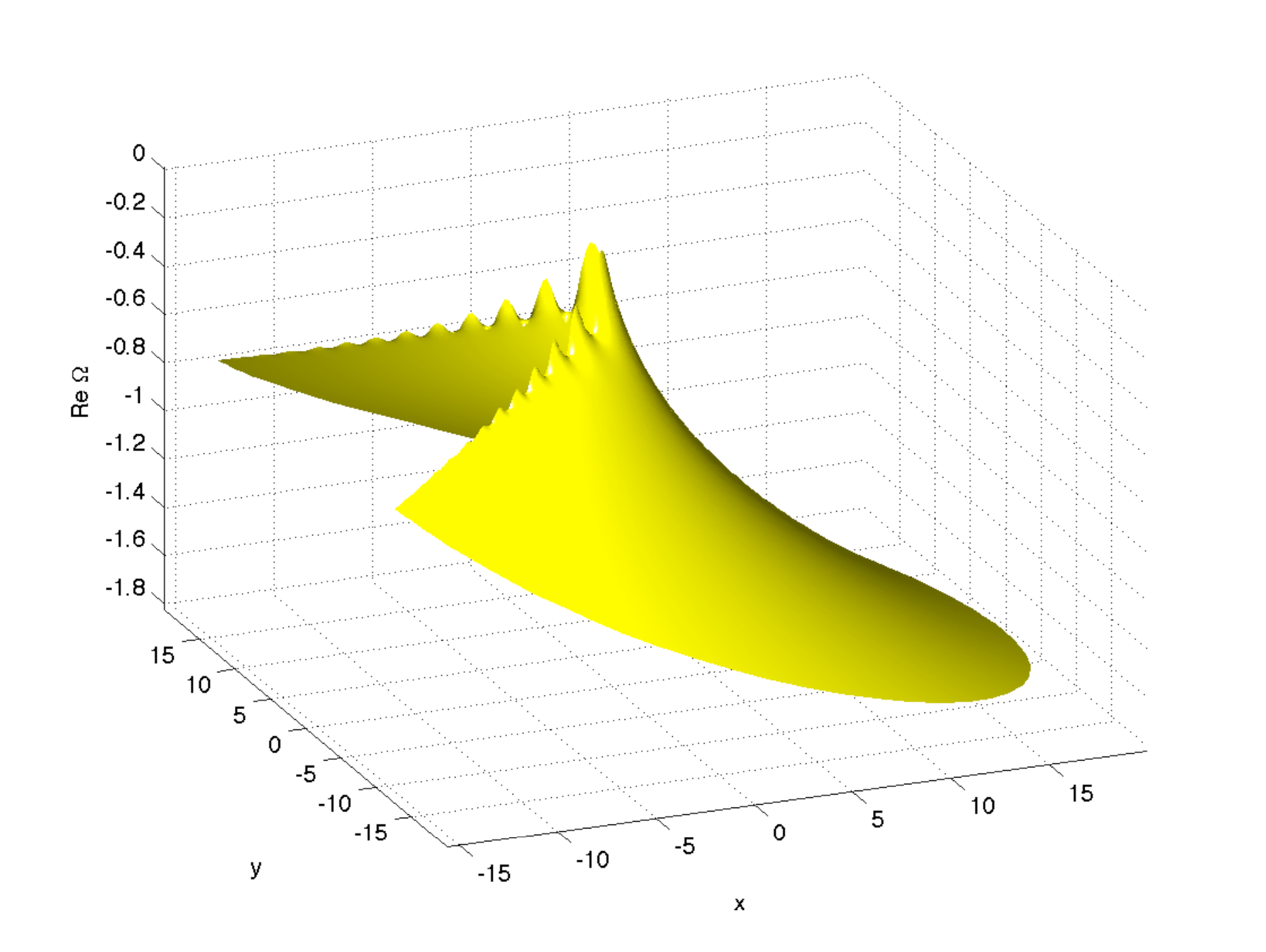}
\caption{Real part of the \emph{tritronqu\'ee} solution in 
the sector $r<20$ and $|\phi|<4\pi/5-0.05$.}
\label{figp1re}
\end{figure}

\begin{figure}[!htb]
\centering
\includegraphics[width=5in]{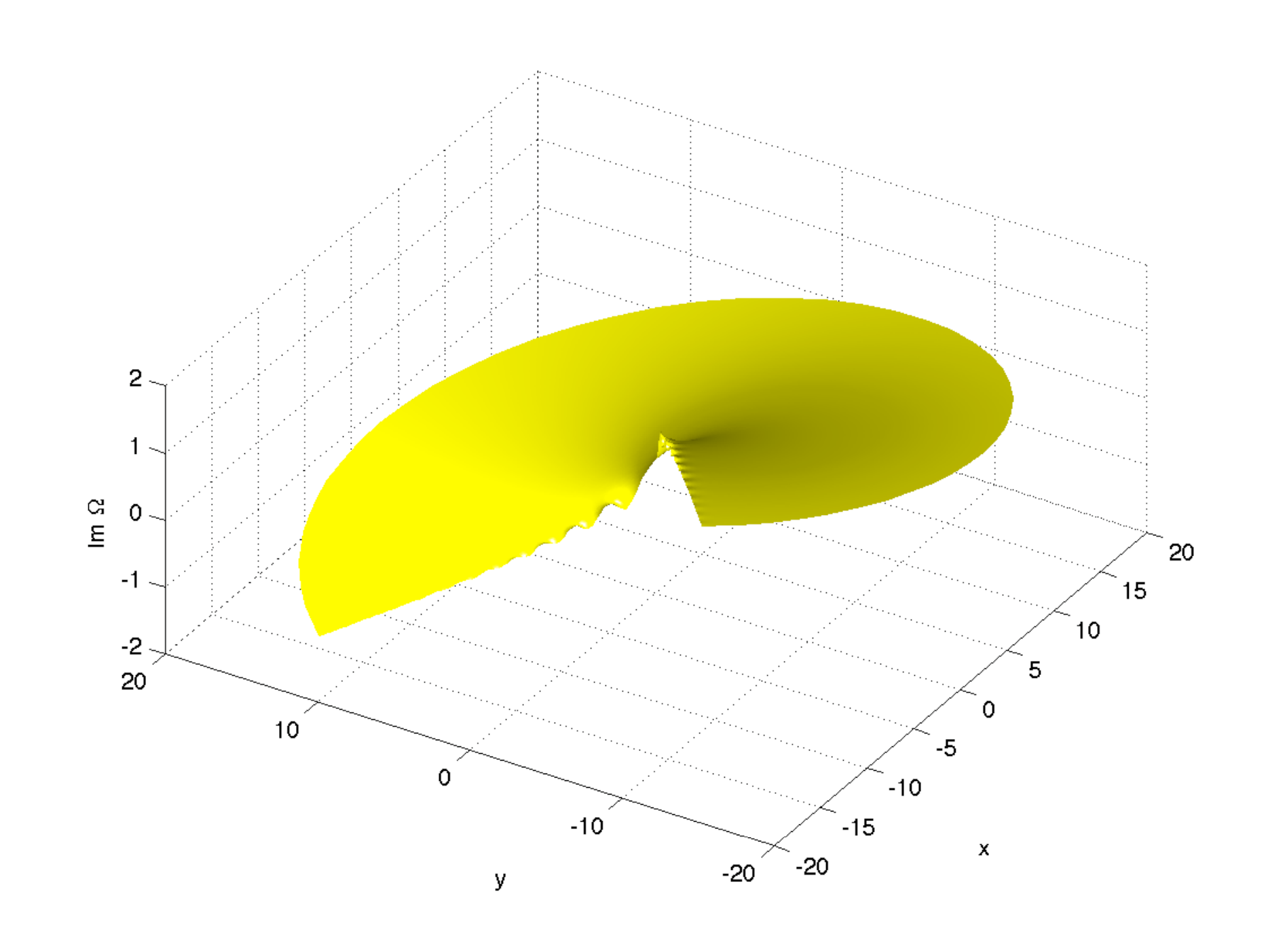}
\caption{Imaginary part of the \emph{tritronqu\'ee} solution in 
the sector $r<20$ and $|\phi|<4\pi/5-0.05$.}
\label{figp1im}
\end{figure}

To solve the resulting boundary value problem for the P-I equation is, however, 
computationally expensive since we have to solve an equation in 2 real 
dimensions
iteratively. Since the solution we want to construct is holomorphic there, we can 
instead solve the harmonicity condition (the two dimensional Laplace 
equation) for the given boundary conditions. To this end we introduce polar 
coordinates $r$, $\phi$ and use a spectral grid as described in 
\cite{trefethen1}: the main point is a doubling of the interval 
$r\in[0,R]$ to $[-R,R]$ to allow for a better distribution of the 
Chebychev collocation points. Since we work with values of 
$\phi<\phi_{0}$, we cannot use the usual Fourier series approach for 
the azimuthal coordinate. Instead we use again a Chebychev 
collocation method. The found solution in the 
considered domain is shown in Fig.~\ref{figp1re} and 
Fig.~\ref{figp1im}.
The quality of the solution can be tested by plugging the found 
solution to the Laplace equation into the P-I equation. Due to the 
low resolution and problems at the boundary, the accuracy is 
considerably lower in the two dimensional case than  on the lines. 
This is, however, not a problem since we need only the one dimensional 
solutions for quantitative comparisons with NLS solutions. The two dimensional solutions 
give nonetheless strong numerical evidence for the conjecture that 
the \emph{tritronqu\'ee} solution has globally no poles in the sector 
$|\arg\zeta|<4\pi/5$. 

\setcounter{equation}{0}
\setcounter{theorem}{0}
\section{Critical behavior in NLS and the  {\it tritronqu\'ee} solution of P-I: numerical results}\label{sec6}\par
In this section we will compare the numerical solution of the 
focusing NLS equation for two examples of initial data for values of 
$\epsilon$ between $0.1$ and $0.025$ with the asymptotic solutions 
discussed in the previous sections, the semiclassical solution 
up to the breakup and the 
tritronqu\'ee solution to the Painlev\'e I equation. The numerical 
approach to solve the NLS equation is discussed in detail in 
\cite{numart1d}. For values of $\epsilon$ below $0.04$ we have to use 
Krasny filtering \cite{krasny} (Fourier coefficients with an absolute 
value below $10^{-13}$ are put equal to zero to avoid the excitation 
of unstable modes). With double precision arithmetic we could thus 
reach $\epsilon=0.025$, but could not go below. 

\subsection{Initial data}
We consider initial data where $u(x,0)$ has a single positive hump, 
and where $v(x,0)$ is monotonously decreasing.
For initial data of the form $u(x,0)=A^2(x)$ and $v(x,0)=0$ where the function
$A(x)$ is analytic with a single positive hump with maximum  value $A_0$, the 
semiclassical solution of NLS follows from (\ref{hodo1}) with $f(u,v)$  given by 
\beq\label{eq1a}
f(u,v)=2\Im\left(\int\limits^{-\frac{1}{2}v+i\sqrt{u}}_{iA_0}d\eta\,\rho(\eta)\sqrt{(\eta+\frac{1}{2}v)^2+u}\right)
\eeq
where
\[
\rho(\eta)=\dfrac{\eta}{\pi}\int_{x^-(\eta)}^{x^+(\eta)}\dfrac{dx}{\sqrt{A^2(x)+\eta^2}},
\]
and where $x^{\pm}(\eta)$ are defined by $A(x^{\pm}(\eta))=i\eta$. 
The formula (\ref{eq1a}) follows from results by S.Kamvissis, K.McLaughlin and 
P.Miller in \cite{kam2}.

From  $f(u,v)$ it is straightforward to recover  the initial data from the equations
\begin{equation}
\label{ID}
x=f_u,\quad f_v=0.
\end{equation}
Numerically we study the critical behavior of two classes of initial data, 
one symmetric with respect to $x$ which were used in \cite{mi}, and 
initial data without symmetry with respect to $x$ which are built 
from the initial data studied in \cite{to1}. For the former class 
the corresponding exact solution of focusing NLS is known in terms of a 
determinant. Nonetheless we integrate the NLS equation for these 
initial data numerically since this approach is not limited to 
special cases, but can be used for general smooth Schwartzian initial 
data as in the latter case. 


\subsubsection{Symmetric  initial data}
We consider the particular class of  initial data data
\begin{equation}
\label{symmdata}
u(x,t=0)=A_{0}^{2}\text{sech}^2\, x,\quad v(x,t=0)=-\mu\,\text{tanh} \, x,\quad \mu\geq 0.
\end{equation}
Introducing the quantity 
\[
M=\sqrt{\frac{\mu^2}{4}-A_{0}^{2}},
\]
we find that  the semiclassical solution for these 
initial data follows from (\ref{hodo1}) with
\begin{equation}
\label{eq11a}
\begin{split}
f(u,v)=&\dfrac{\mu}{2}v-\dfrac{1}{4}(v-2M)\Delta_+
-\dfrac{1}{4}(v+2M)\Delta_--\dfrac{1}{2}u\log u\\
&+\dfrac{1}{2}u\log\left[(-\frac{1}{2}v+M+\Delta_+)(-\frac{1}{2}v-M+\Delta_-)\right]
\end{split}
\end{equation} 
where 
\[
\Delta_{\pm}=\left((-\dfrac{1}{2}v\pm M)^2+u\right)^{\frac{1}{2}}.
\]
For $\mu=0$ we recover  the Satsuma-Yajima \cite{saya} initial data that were 
studied numerically in \cite{mi}. The function $f(u,v)$ takes the form  
\beq\label{eq1}
f(u,v)=\Re \left[ \left( - \frac{v}{2} +i\, A_0\right) 
\sqrt{u+\left(- \frac{v}{2} +i\, A_0\right)^2} + u \, \log\frac{ 
(- \frac{v}{2} +i\, 
A_0)^{2}+\sqrt{(-\frac{v}{2}+iA_{0})^{2}+u}}{\sqrt{u}}\right],
\eeq
which can also be recovered from (\ref{eq1a}) by setting $\rho=i$.
The critical point is given by
\begin{equation}
    u_{0}=2A_0^{2},\quad v_{0}=0,\quad x_{0}=0,\quad 
    t_{0}=\frac{1}{2A_0}.
    \label{eq3}
\end{equation}
Furthermore we have
\begin{equation}
    f_{uuu}^{0}=0,\quad f_{uuv}^{0}=\frac{1}{4A_0^{3}},\quad r =  
    4A_0^{3},\quad \psi=0,
    \label{eq4}
\end{equation}
where $r$, $\psi$ are defined in (\ref{rpsi0}). 
For $A_{0}=1$ the initial data (\ref{symmdata}) coincides with the one 
studied in  \cite{to1} by
A.Tovbis, S.Venakides and X.Zhou.
In the particular case $\mu=2$, $A_{0}=1$ the function $f(u,v)$ in 
(\ref{eq1a}) simplifies to
\begin{equation}
\label{eq11a2}
f(u,v)= v - \frac{v}2\, \sqrt{\frac{1}{4}v^2+u} +u\log
\left[\dfrac{-\frac{1}{2}v+\sqrt{\frac{1}{4}v^2+u} }
{\sqrt{u}}\right].
\end{equation}
In this case the critical point is given by
\[
v_{0}=0,\;\;u_{0}=2+\mu,
\quad t_{0}=\dfrac{1}{2+\mu},\;\;x_{0}=0.
\]
Furthermore, 
\[
f_{uuu}^{0}=0,\quad f_{uuv}^{0}=\dfrac{2}{(\mu+2)^3},
\quad
 r = \dfrac{(\mu+2)^3} {2},\quad \psi=0.
\]
\subsubsection{Non-symmetric initial data}
Recall that we are 
interested here in Cauchy data in the Schwartz class of rapidly 
decreasing functions.
The above initial data are symmetric with respect to $x$, $u$ is an even 
and $v$ an odd function in $x$.  To obtain a situation which is manifestly not 
symmetric, we use the fact that if $f$ is a solution to (\ref{eq-f}),
the same holds 
for derivatives and anti-derivatives of $f$ with respect to $v$ and
for any linear combination of those. If 
$f_{v}$ is an even function in $v$, this will obviously not be the 
case for a linear combination of $f$ and $f_{v}$. 

As a specific example, we consider the linear combination 
\[
f =  f_{1}+\alpha f_{2},\quad \alpha=const,
\]
 where $f_1$ coincides with  (\ref{eq11a2})
and
\begin{equation}
 \label{f1}
    f_{2}  =2\,u\,\sqrt{\frac{1}{4}v^2+u} - \frac{2}{3}\left(
    \frac{1}{4}v^2+u\right)^{3/2} +u\, v\log
\left[\dfrac{-\frac{1}{2}v+\sqrt{\frac{1}{4}v^2+u} }
{\sqrt{u}}\right]. 
\end{equation}
The function $f_{2}$ is obtained from the integration
$f_{2,v}=f_{1}-v$.
The critical point is given in this case by
$$
u_{0}=4(1-16\alpha^{2}),\quad v_{0}=-16\alpha,\quad 
x_{0}=\frac{1}{2}\log \frac{1+4\alpha}{1-4\alpha},\quad 
t_{0}=\frac{1}{4}-\frac{\alpha}{2}\log \frac{1+4\alpha}{1-4\alpha};$$ 
thus we have $|\alpha|<1/4$. Furthermore, 
\[
f_{uuv}^{0}=-\frac{4\alpha^{2}-1/8}{4\sqrt{1-16\alpha^{2}}},\quad 
f_{uuu}^{0}=\frac{\alpha}{4\sqrt{1-16\alpha^{2}}},
\]
such that
\[
r=8u_0,\quad \psi=-\arctan\dfrac{\alpha\sqrt{1-16\alpha^2}}{1/8-4\alpha^2}.
\]

\begin{figure}[!htb]
\centering
\includegraphics[width=5in]{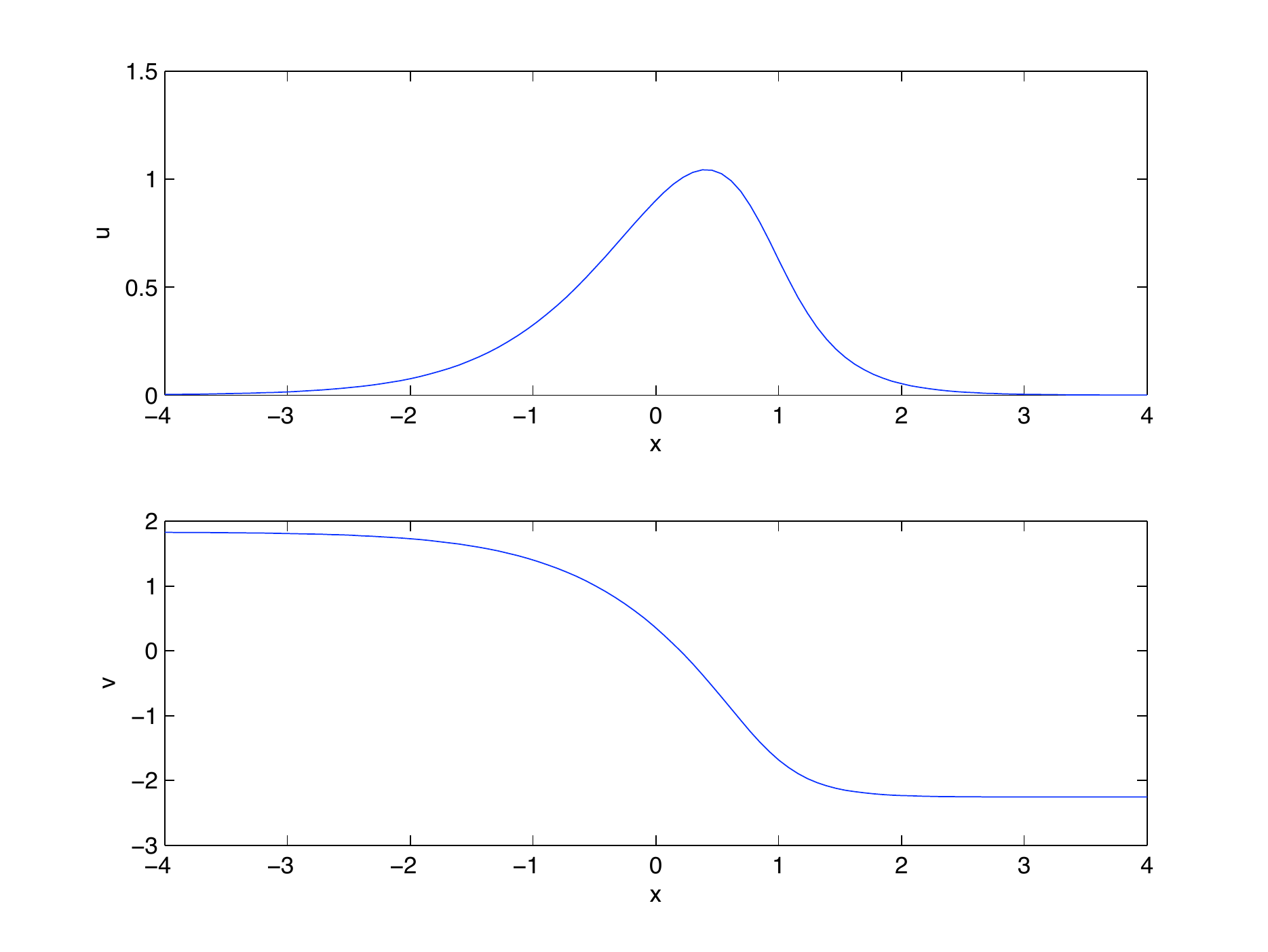}
\caption{Initial data for the NLS equations without symmetry with 
respect to $x$.}
\label{fignlsasymini}
\end{figure}

We determine the initial data corresponding to $f$ for a 
given value of $|\alpha|<1/4$ by solving (\ref{ID}) for 
$u$, $v$ in dependence of $x$. This is 
done numerically by using the algorithm of \cite{optim} which is 
implemented as the function \emph{fminsearch} in Matlab. The 
algorithm provides an iterative approach which converges in our case 
rapidly if the starting values are close enough to the solution, 
which is achieved by choosing $u$ and $v$ corresponding to 
$f_{0}$ as an initial guess. 
For $\alpha$ close to $1/4$ we observe 
numerically a steepening of the initial pulse which will lead to a shock 
front in the limit $\alpha\to 1/4$. For the computations presented 
here, we consider the case $\alpha=0.1$ which leads to the initial 
data shown in Fig.~\ref{fignlsasymini}.

The initial data are computed in the way described above to the 
order of the Krasny filter on the interval $x\in [-15,11]$ on 
Chebychev collocation points. Standard interpolation via Chebychev 
polynomials is then used to interpolate the resulting data to a 
Fourier grid. To avoid a Gibbs phenomenon at the interval ends due 
to the non-periodicity of the data, we 
use a Fourier grid on the interval $[-10\pi,10\pi]$ to ensure that 
the function $u$ takes values of the order of the Krasny filter. 
For $x<-15$ and 
$x>11$, the function $u$ is exponentially small which implies the 
zero-finding algorithm will no longer provide the needed precision. 
Thus we determine the exponential tails of the solution to leading 
order analytically. We find for $x\to-\infty$
\begin{align}
    u &= v_{+}^2\exp\left(\frac{2(x-\alpha v_{+})}{\alpha v_{+}+1}
    \right),\nonumber\\ 
    v &= v_{+} - v_{+}\exp\left(\frac{2(x-\alpha v_{+})}{\alpha v_{+}+1}
    \right)
    \left(2\log(v_{+})+\frac{2(x-\alpha v_{+})}{\alpha v_{+}+1}\right)
    \label{infmin},
\end{align}
and for $x\to+\infty$
\begin{align}
    u &= v_{-}^2\exp\left(-\frac{2(x+\alpha)}{\alpha v_{-}+1}\right)\nonumber\\
    v &= v_{-} - v_{-}\exp\left(-\frac{2(x+\alpha)}{\alpha v_{-}+1}\right)
    \left(2\log(-v_{-})-\frac{2(x+\alpha)}{\alpha v_{-}+1}\right)
    \label{infpl},
\end{align}
where $v_{\pm}=(\sqrt{1\pm \alpha}-1)/\alpha$.
The initial data for the NLS equation in the form 
$\Psi=\sqrt{u}\exp(iS/\epsilon)$ are then found by integrating $v$ 
on the Chebychev grid by standard integration of Chebychev 
polynomials. The exponential tails for $S$ 
follow from (\ref{infmin}) and 
(\ref{infpl}). The matching of the tails to the Chebychev interpolant 
is not smooth and leads to a small Gibbs phenomenon. The Fourier 
coefficients decrease, however, to the order of the Krasny filter which is 
sufficient for our purposes. Thus we obtain the non-symmetric initial 
data with roughly the same precision as the analytic symmetric data.

\subsection{Semiclassical solution}
For times $t\ll t_{c}$, the semiclassical solution gives a very 
accurate asymptotic description for the NLS solution. The situation 
is similar to the Hopf and the KdV equation \cite{gk1}. We find for the 
symmetric initial data for $t=t_{c}/2$
that the 
$L_{\infty}$ norm of the difference between the solutions decreases as 
$\epsilon^{2}$. More precisely a linear regression analysis in the 
case of symmetric initial data (for the 
values $\epsilon=0.03,0.04,\ldots,0.1$) for the 
logarithm of this norm leads to an error proportional to 
$\epsilon^{a}$ with $a=1.94$, a correlation coefficient 
$r=0.9995$ and standard error $\sigma_{a}=0.03$. In the non-symmetric 
case, we find $a=1.98$, $r=0.999996$ and $\sigma_{a}=0.003$.

\begin{figure}[!htb]
\centering
\includegraphics[width=5in]{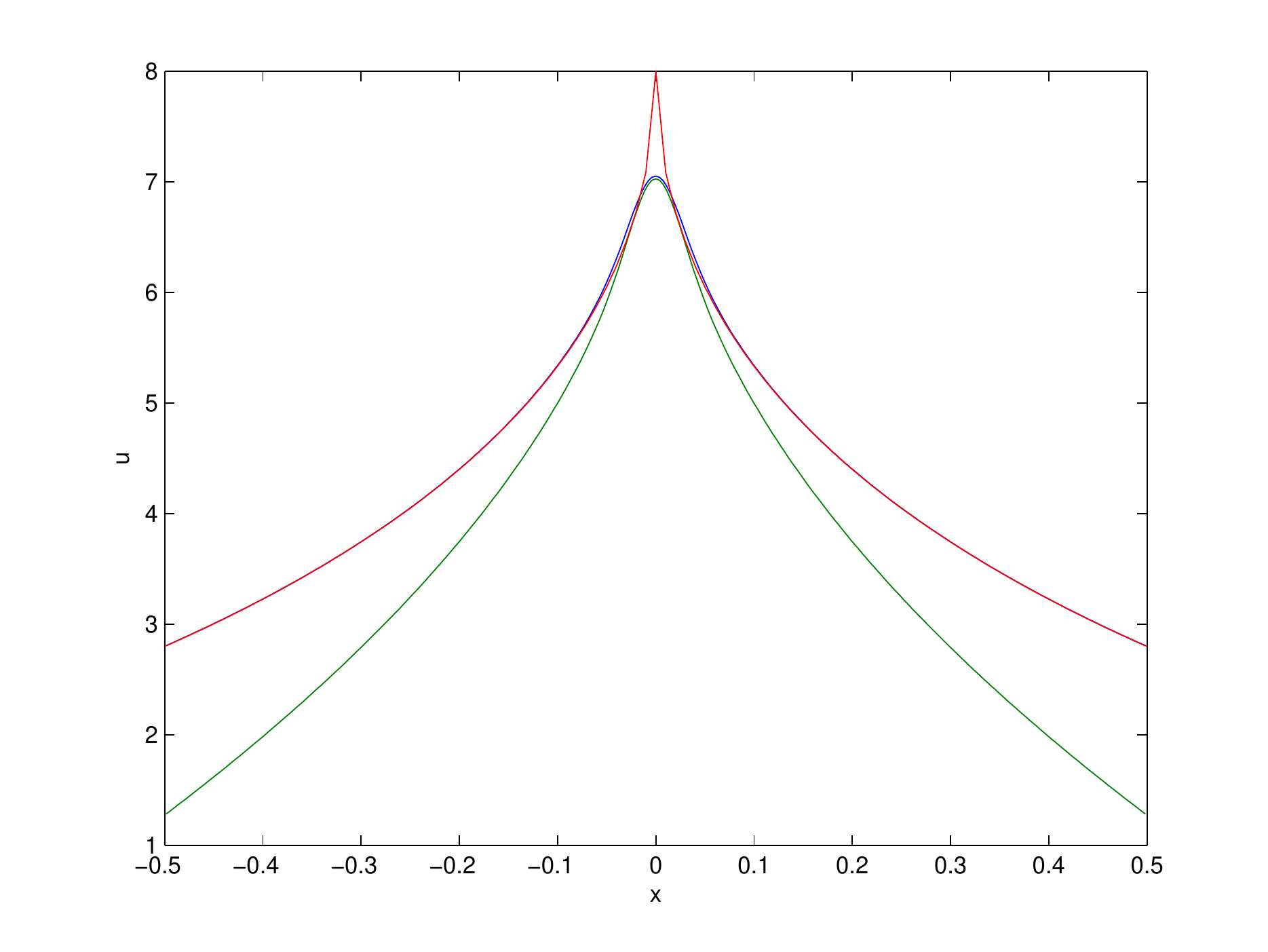}
\caption{The blue line is the function $u$ of the 
solution to the focusing NLS equation for the 
initial data $u(x,0)=2\,\mbox{sech}\,x$ and $\epsilon=0.04$ at the 
critical time, 
and the red line is the corresponding  semiclassical solution
given by formulas (\ref{hodo}). The green line gives 
the multiscales solution via the tritronqu\'ee solution of the 
Painlev\'e I equation.}
\label{fignlscritical}
\end{figure}

\begin{figure}[!htb]
\centering
\includegraphics[width=5in]{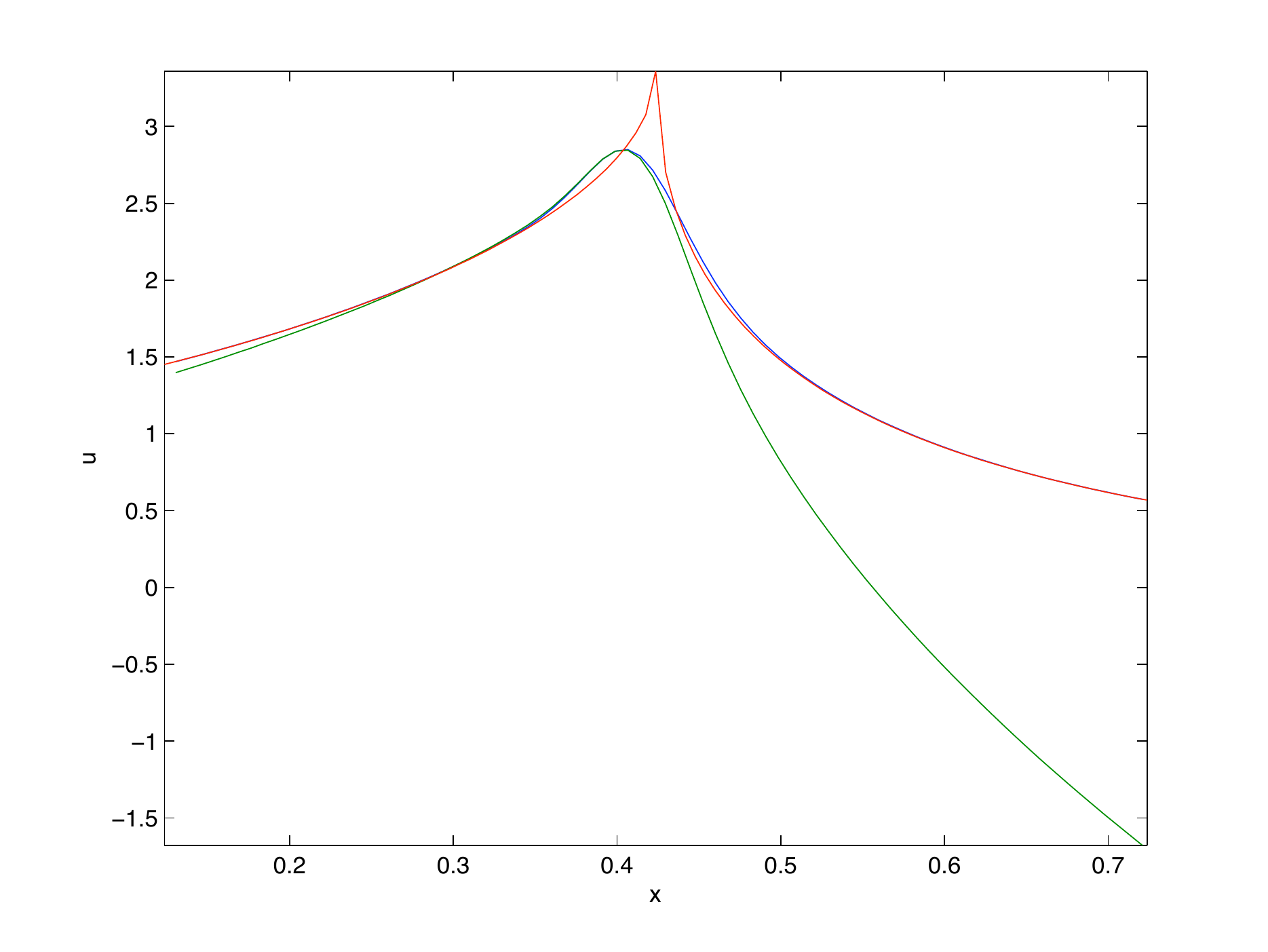}
\caption{The blue line is the function $u$ of the 
solution to the focusing NLS equation for the non-symmetric
initial data and $\epsilon=0.04$ at the 
critical time, 
and the red line is the corresponding  semiclassical solution
given by formulas (\ref{hodo}). The green line gives 
the multiscales solution via the tritronqu\'ee solution of the 
Painlev\'e I equation.}
\label{fignlscriticalasym}
\end{figure}

Close to the critical time the semiclassical solution only provides a 
satisfactory description of the NLS solution for large values of 
$|x-x_{c}|$. In the breakup region it fails to be accurate since 
it develops a cusp at $x_{c}$ whereas the NLS solution stays smooth. 
This behavior can be well seen in Fig.~\ref{fignlscritical} for the symmetric initial 
data.
The largest difference between the semiclassical and the NLS solution 
is always at the critical point. We find that the $L_{\infty}$ norm 
of the difference scales roughly as $\epsilon^{2/5}$ as suggested by 
the Main Conjecture. More precisely we find a scaling proportional to 
$\epsilon^{a}$ with $a =0.38$ and $r=0.999997$ and 
$\sigma_{a}=4.2*10^{-4}$. For the non-symmetric initial data, we find 
$a=0.36$, $r=0.9999$ and $\sigma_{a}=0.002$.
The corresponding plot for $u$ can be seen in  
Fig.~\ref{fignlscriticalasym}.

\begin{figure}[!htb]
\centering
\includegraphics[width=5in]{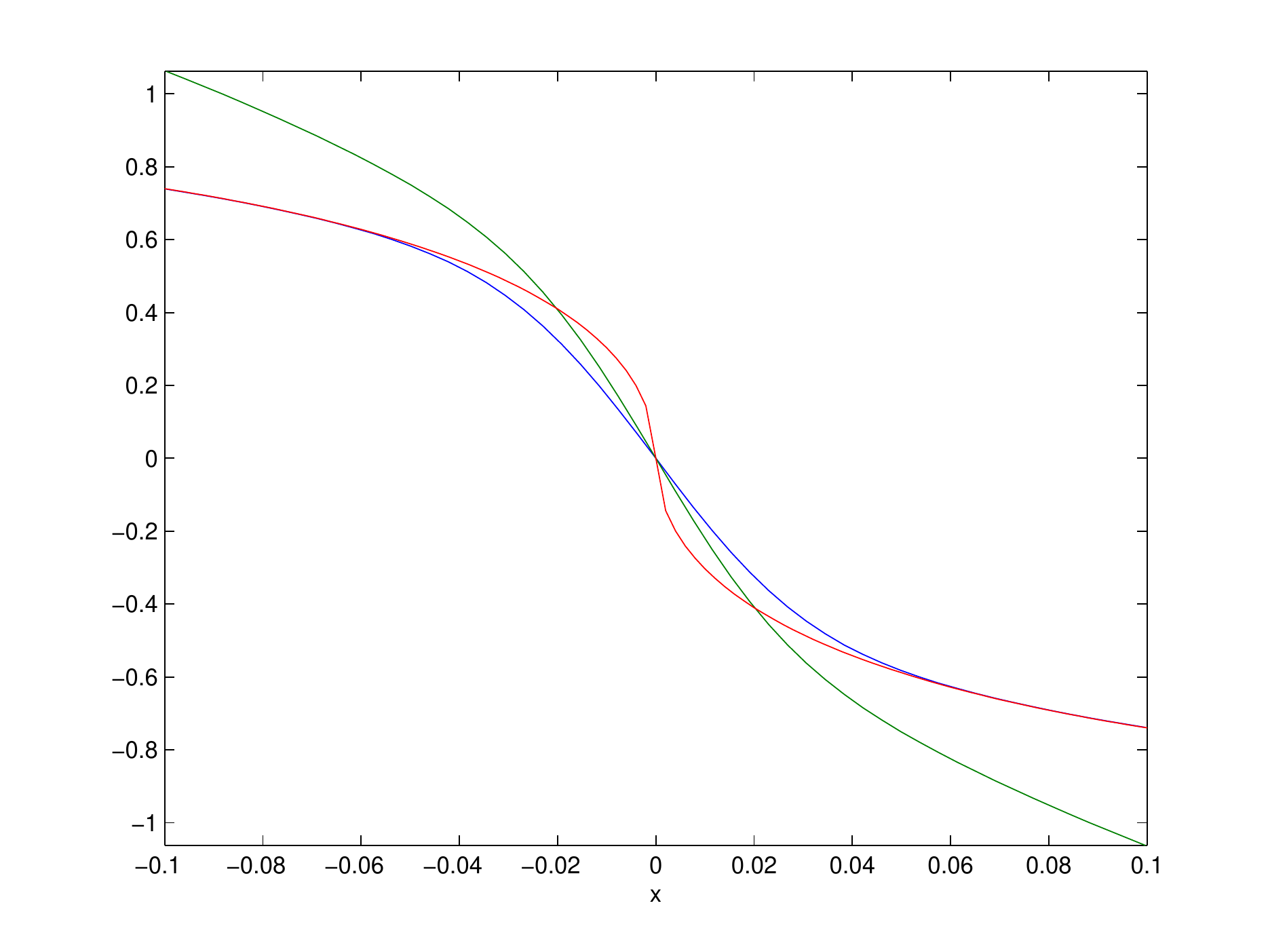}
\caption{The blue line is the function $v$ of the 
solution to the focusing NLS equation for the 
initial data $u_0(x)=2\,\mbox{sech}\,x$ and $\epsilon=0.04$ at the 
critical time, 
and the red line is the corresponding  semiclassical solution
given by formulas (\ref{hodo}). The green line gives 
the multiscales solution via the tritronqu\'ee solution of the 
Painlev\'e I equation.}
\label{fignlscriticalv}
\end{figure}

\begin{figure}[!htb]
\centering
\includegraphics[width=5in]{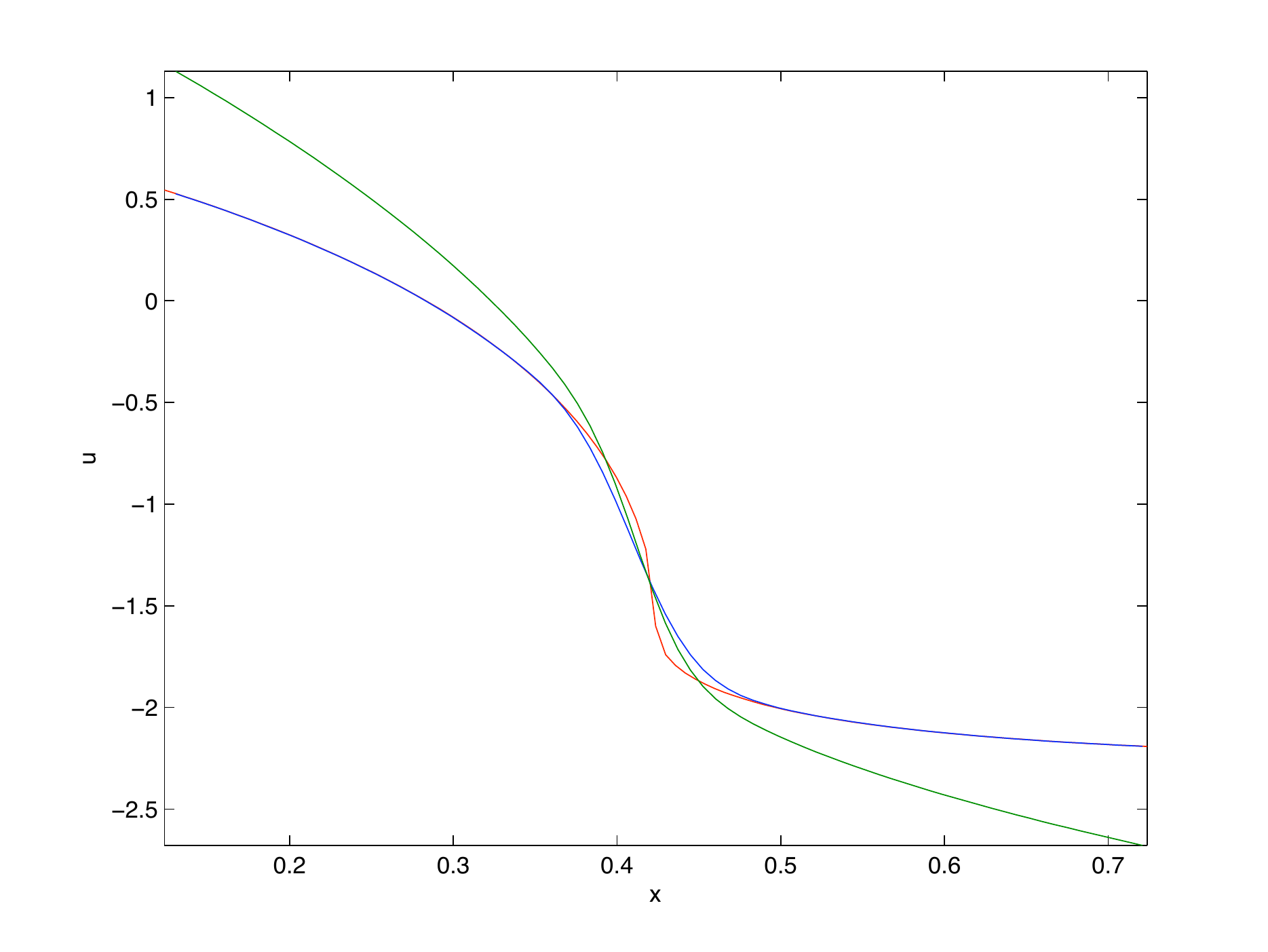}
\caption{The blue line is the function $v$ of the 
solution to the focusing NLS equation for the non-symmetric
initial data and $\epsilon=0.04$ at the 
critical time, 
and the red line is the corresponding  semiclassical solution
given by formulas (\ref{hodo}). The green line gives 
the multiscales solution via the tritronqu\'ee solution of the 
Painlev\'e I equation.}
\label{fignlscriticalasymv}
\end{figure}

The function $v$ for the same situation as in 
Fig.~\ref{fignlscritical} is shown in Fig.~\ref{fignlscriticalv}. It 
can be seen that the semiclassical solution is again a satisfactory 
description for $|x-x_{c}|$ large, but fails to be accurate close to 
the breakup point. 
The phase for the non-symmetric initial data can be seen in 
Fig.~\ref{fignlscriticalasymv}. 
In the following we will always study the scaling 
for the function $u$ without further notice. 

\subsection{Multiscales solution}
\begin{figure}[!htb]
\centering
\includegraphics[width=5in]{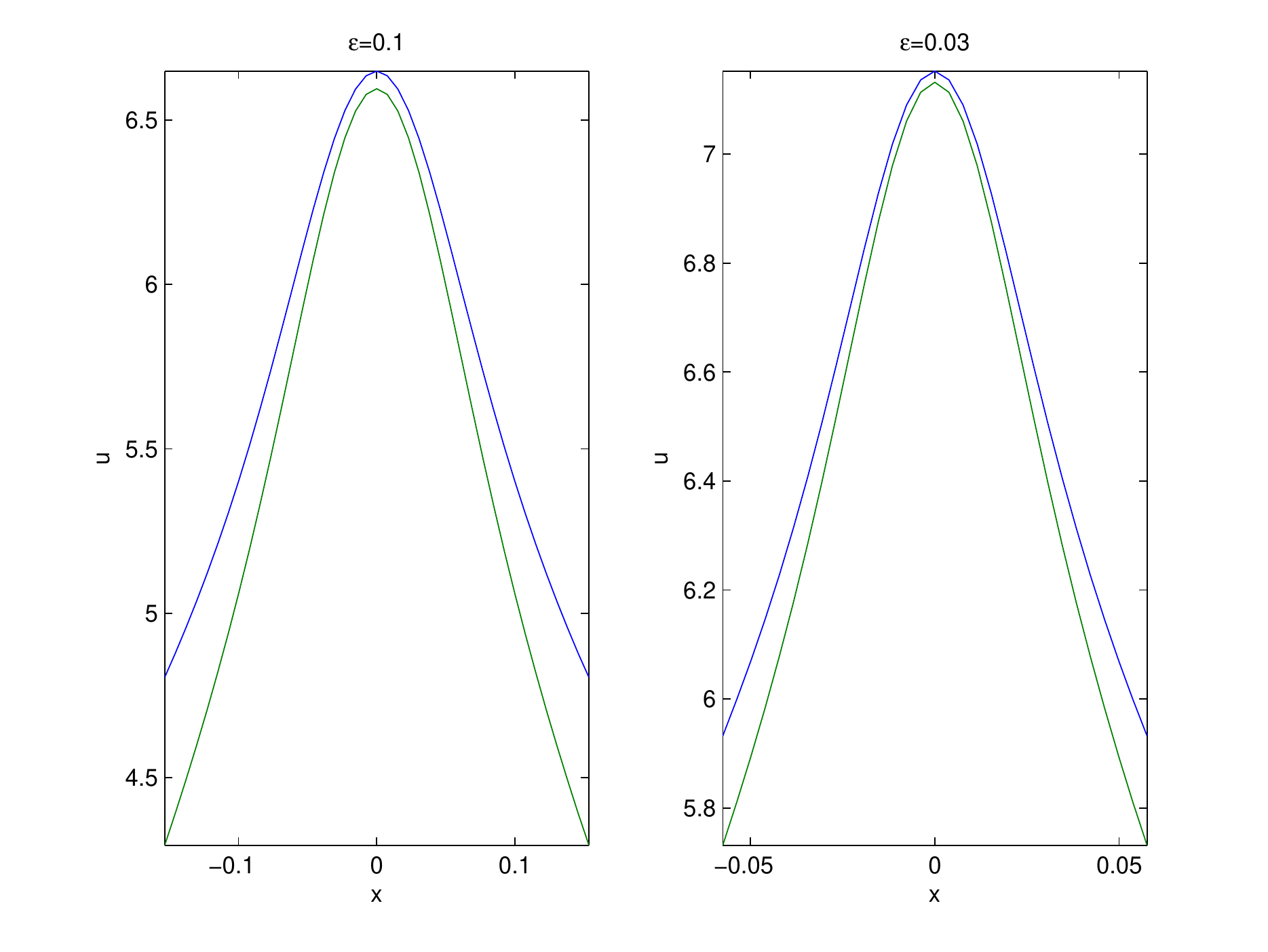}
\caption{The blue line is the solution to the focusing NLS equation for the 
initial data $u_0(x)=2\,\mbox{sech}\,x$ at the 
critical time, and the green line gives 
the multiscales solution via the tritronqu\'ee solution of the 
Painlev\'e I equation. The plots are shown for two values of 
$\epsilon$ at the critical time.}
\label{fignlsc2e}
\end{figure}

It can be seen in Fig.~\ref{fignlscritical} and 
Fig.~\ref{fignlscriticalv} 
that the multiscales solution (\ref{main}) in terms 
of the tritronqu\'ee solution to the Painlev\'e I equation gives a 
much better asymptotic description to the NLS solution at breakup 
close to the breakup point than the semiclassical solution for the 
symmetric initial data. For 
larger values of $|x-x_{c}|$, the semiclassical solution provides, 
however, the better approximation. The rescaling of the coordinates in 
(\ref{main}) suggests to consider the difference between the NLS 
and the multiscales solution in an interval $[-\gamma 
\epsilon^{4/5},\gamma \epsilon^{4/5}]$ (we choose here $\gamma=1$, 
but within numerical accuracy the result does not depend on varying $\gamma$ around this 
value). These intervals can be seen in 
Fig.~\ref{fignlsc2e}. We find that the $L_{\infty}$ norm 
of the difference between these solutions in this interval scales 
roughly like $\epsilon^{4/5}$. More precisely we have a scaling 
$\epsilon^{a}$ with $a = 0.76$ ($r=0.998$ and $\sigma_{a}=0.019$).

\begin{figure}[!htb]
\centering
\includegraphics[width=5in]{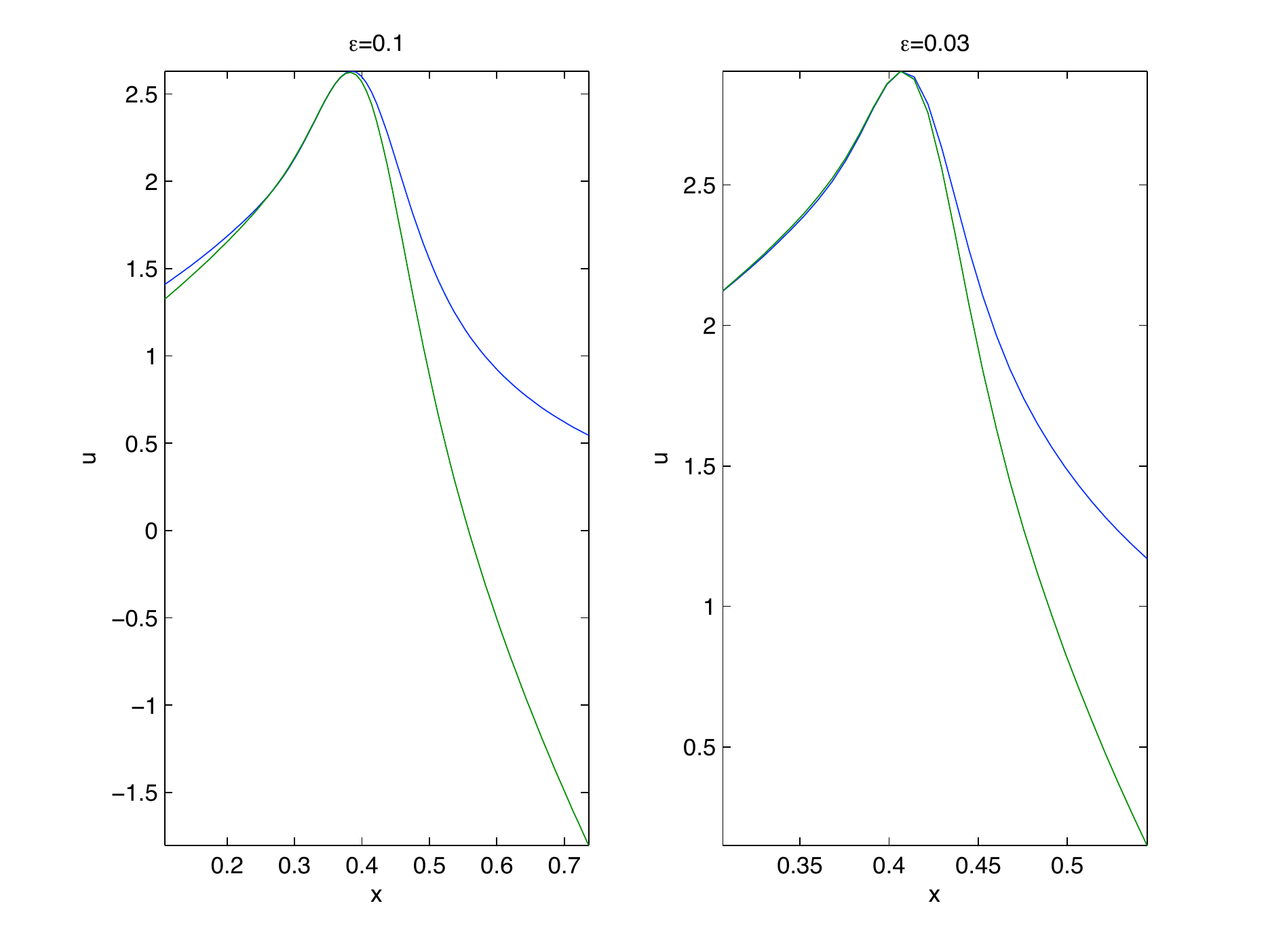}
\caption{The blue line is the solution to the focusing NLS equation for 
the non-symmetric initial data at the 
critical time, and the green line gives 
the multiscales solution via the tritronqu\'ee solution of the 
Painlev\'e I equation. The plots are shown for two values of 
$\epsilon$ at the critical time.}
\label{fignlsc2easym}
\end{figure}

For the non-symmetric initial data, the situation at the critical 
point can be seen in Fig.~\ref{fignlscriticalasym} and 
Fig.~\ref{fignlscriticalasymv}. Again the multiscales solution 
(\ref{main}) gives a much better description close to the critical 
point than the semiclassical solution. However, the approximation is 
here much better on the side with weak slope for $u$ than on the side with 
strong slope. We consider again the $L_{\infty}$-norm of the 
difference between the multiscales and the NLS solution in the 
interval $[-\gamma \epsilon^{4/5},\gamma \epsilon^{4/5}]$. The  
scaling behavior of the solution can be seen in 
Fig.~\ref{fignlsc2easym}. For 
$\gamma=1$ we find $a=0.71$, $r=0.998$ and $\sigma_{a}=0.02$. These 
values do not change much for larger $\gamma$. For smaller $\gamma$ 
there are not enough points to provide a valid statistics. The value 
of $a$ smaller than the predicted $4/5$ is seemingly due to the 
strong asymmetry in the quality of the approximation of NLS by the 
multiscales solution as can be seen from 
Fig.~\ref{fignlscriticalasym}. In the considered interval, the 
deviation is already so big that the scaling no longer holds as in 
the symmetric case. To study the scaling with a reliable statistics 
would, however, require the use of a considerably higher resolution 
which would be computationally too expensive.

\begin{figure}[!htb]
\centering
\includegraphics[width=6in]{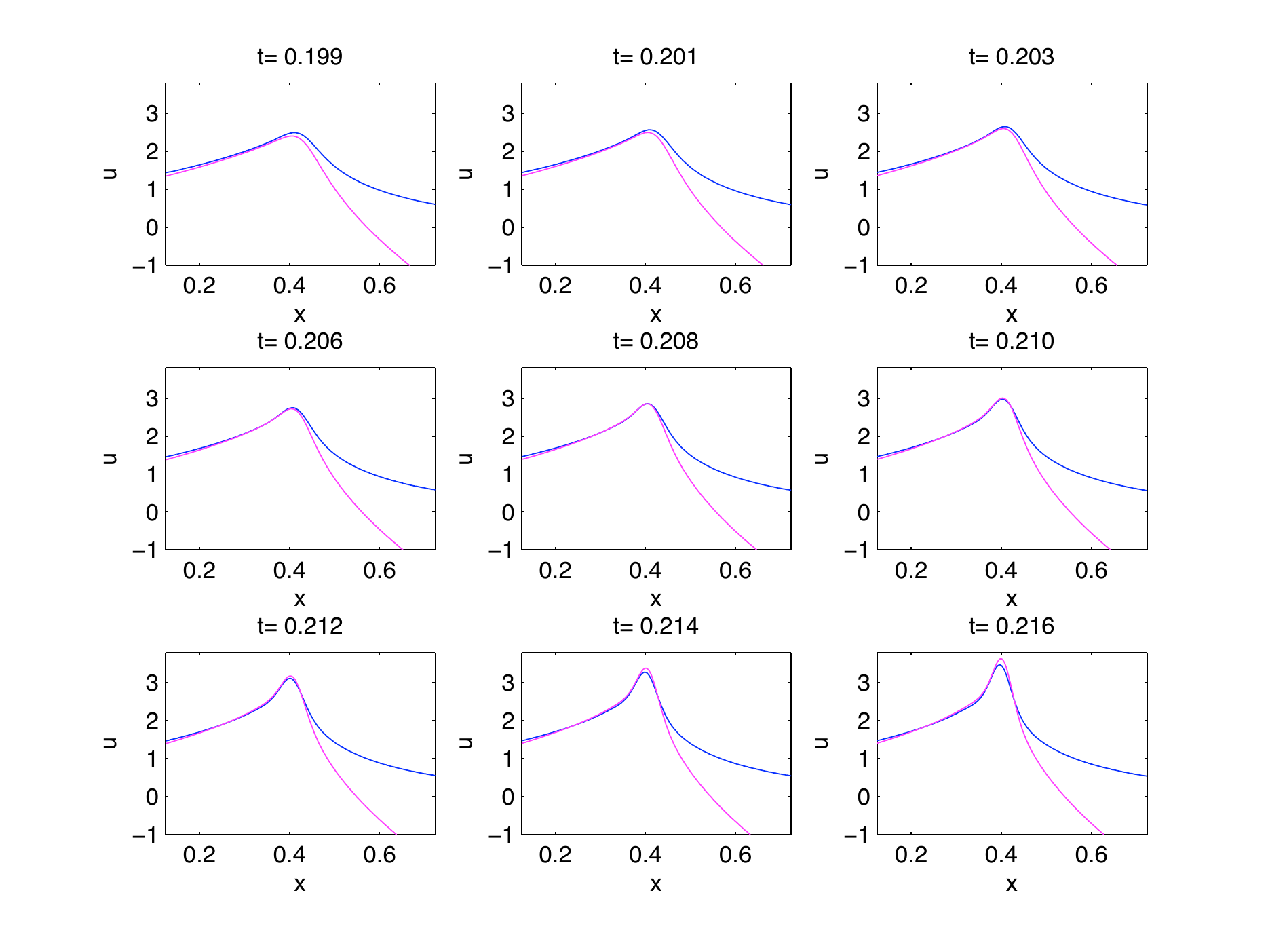}
\caption{The blue line is the solution to the focusing NLS equation for 
the non-symmetric initial data for $\epsilon=0.04$
for various times, and the magenta line gives 
the multiscales solution via the tritronqu\'ee solution of the 
Painlev\'e I equation. The plot in the middle shows the behavior at 
the critical time.}
\label{fignlsasymbreaku_.04}
\end{figure}

Going beyond the critical time, one finds that the real part of the 
NLS solution continues to grow before the central hump breaks up 
into several humps. Notice that the multiscales solution always leads to a function $u$ 
that is smaller than the corresponding function of the NLS solution 
at breakup and before. 
This changes for times after the breakup as can be inferred from 
Fig.~\ref{fignlsasymbreaku_.04} which shows the time dependence of 
the NLS and the corresponding multiscales solution for the 
non-symmetric initial data. The approximation is always best at 
the critical time.

\begin{figure}[!htb]
\centering
\includegraphics[width=5in]{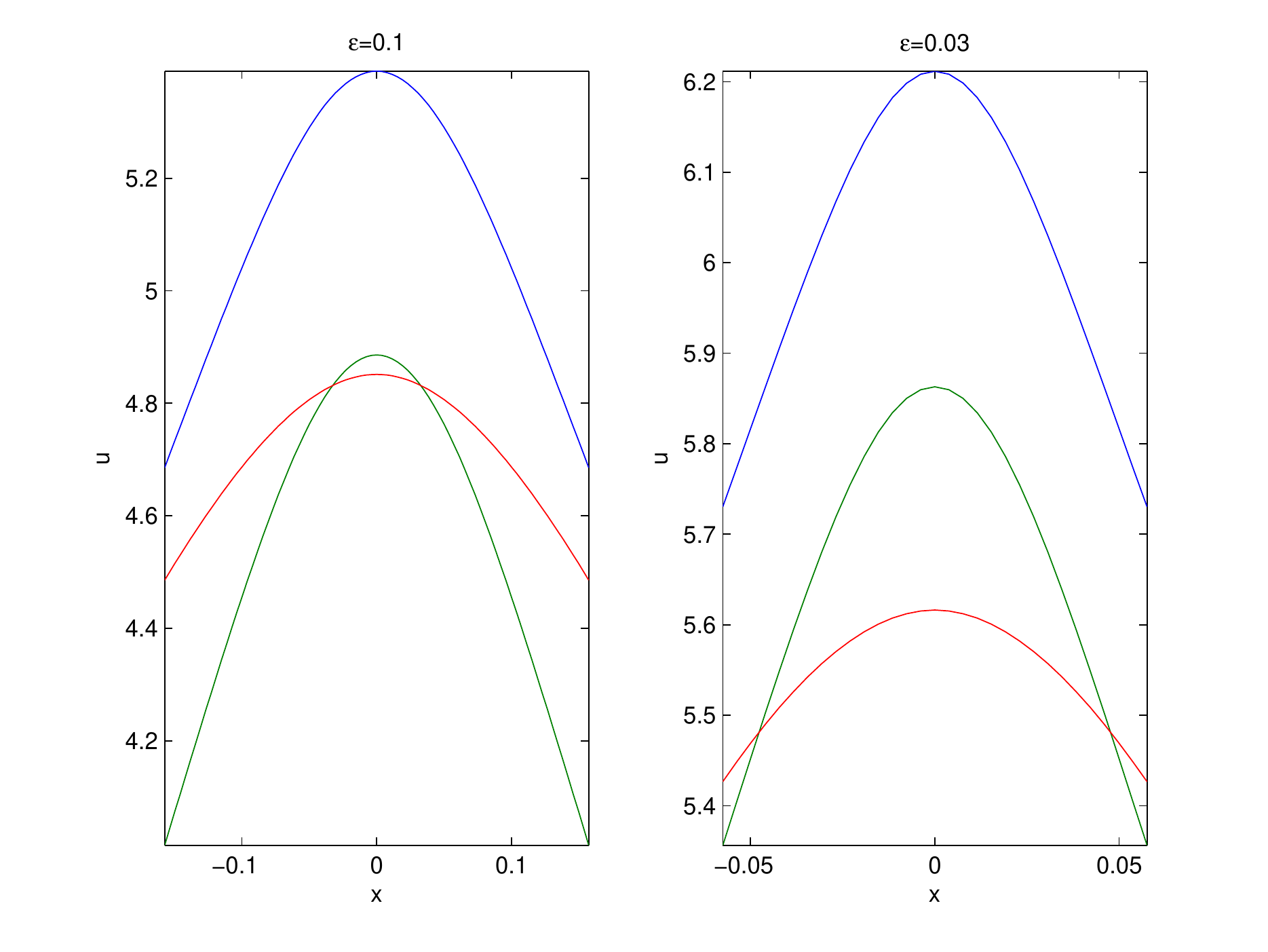}
\caption{The blue line is the solution to the focusing NLS equation for the 
initial data $u_0(x)=2\,\mbox{sech}\,x$, 
and the red line is the corresponding  semiclassical solution
given by formulas  (\ref{hodo}). The green line gives 
the multiscales solution via the tritronqu\'ee solution of the 
Painlev\'e I equation. The plots are shown for two values of 
$\epsilon$ at the corresponding times $t_{-}(\epsilon)$.}
\label{fignlscb2e0}
\end{figure}

To study the quality of the approximation 
(\ref{main}), we use rescaled times. 
The scaling of the coordinates in (\ref{main}) suggests to consider 
the NLS solution close to breakup at the times $t_{\pm}(\epsilon)$ with 
\begin{equation}
    t_{\pm}(\epsilon) = t_{c} +u_{0}/r -\sqrt{(u_{0}/r)^2\pm \epsilon^{4/5}\beta}
    \label{tpm},
\end{equation}
where $\beta$ is a constant (we consider $\beta=0.1$). We will only 
study the symmetric initial data in this context.  Before breakup 
we obtain  the situation shown in Fig.~\ref{fignlscb2e0}. It can be 
seen that the multiscales solution always provides a better 
description close to $x_{c}$ than the semiclassical solution, and 
that the quality improves in this respect with decreasing $\epsilon$. 
We find that the $L_{\infty}$ norm of the difference scales in this 
case as $\epsilon^{a}$ with $a=0.55$ ($r=0.994$ and 
$\sigma_{a}=0.03$).  

\begin{figure}[!t]
\centering
\includegraphics[width=5in]{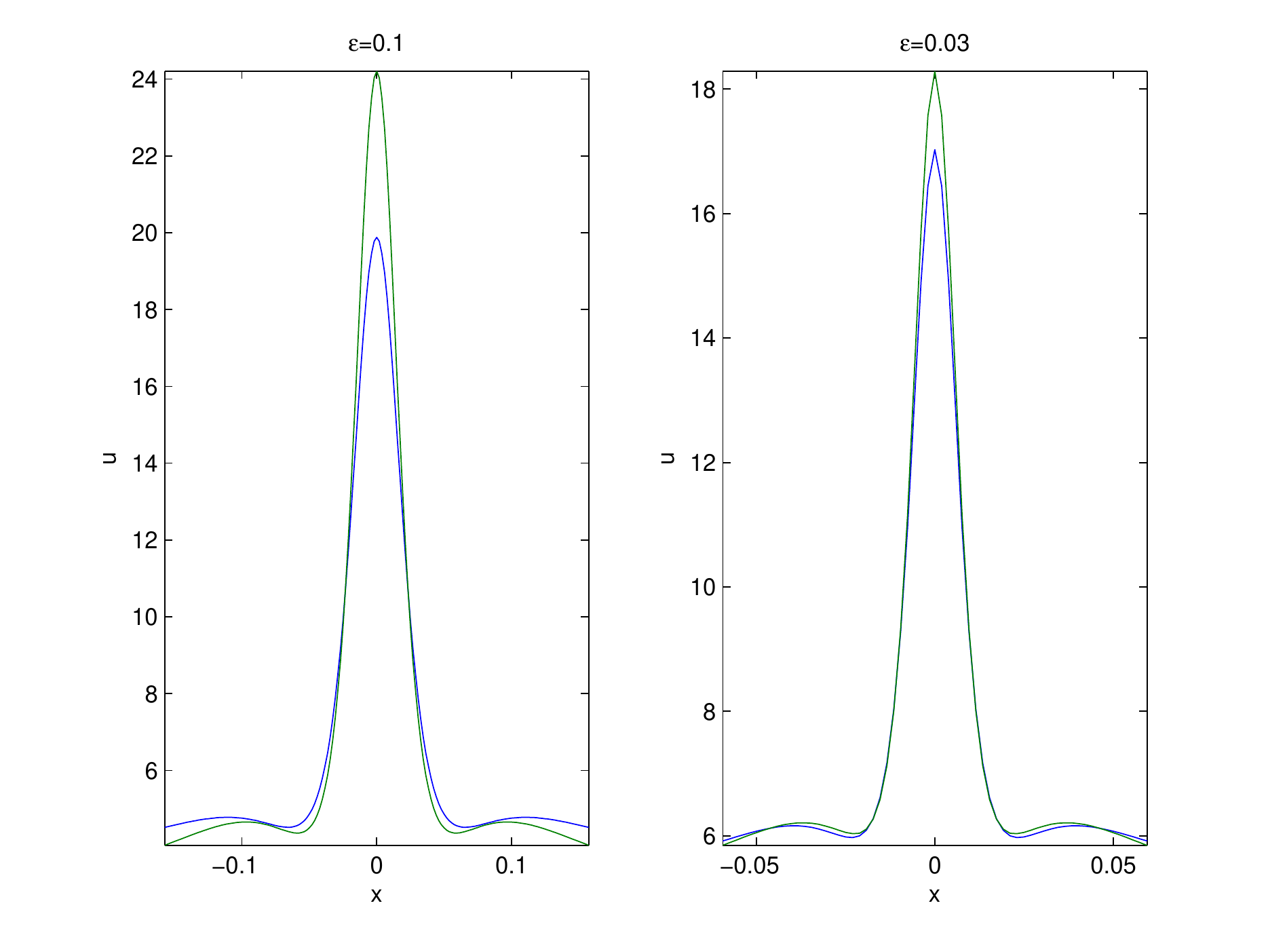}
\caption{The blue line is the solution to the focusing NLS equation for the 
initial data $u_0(x)=2\,\mbox{sech}\,x$, and the green line gives 
the multiscales solution via the tritronqu\'ee solution of the 
Painlev\'e I equation. The plots are shown for two values of 
$\epsilon$ at the corresponding times $t_{+}(\epsilon)$.}
\label{fignlsca2e}
\end{figure}

The situation for times after breakup can be inferred from 
Fig.~\ref{fignlsca2e}. Close to the central region the multiscales 
solution shows a clear difference to the NLS solution. But it is 
interesting to note that the ripples next to the central hump are well 
approximated by the Painlev\'e I solution. The $L_{\infty}$ norm of 
the difference between the two solutions scales roughly like 
$\epsilon$. More precisely we find a scaling $\epsilon^{a}$ with 
$a=1.02$ ($r=0.9999$ and $\sigma_{a}=7.7*10^{-3}$).

\setcounter{equation}{0}
\setcounter{theorem}{0}
\section{Concluding remarks}\label{sec7}\par

In this paper we have started the study of the critical behavior
of generic solutions of the focusing nonlinear Schr\"odinger 
equation. We have formulated the conjectural analytic description of 
this behavior in terms of the {\it tritronqu\'ee} solution to the 
Painlev\'e-I equation restricted to certain lines in the complex plane. We provided analytical as well as numerical evidence supporting our conjecture. In subsequent publications we plan to further study the Main Conjecture of the present paper by applying techniques based, first of all, on the Riemann - Hilbert problem method \cite{kam2, to1, to2} and the theory of Whitham equations (see \cite{gk1} for the numerical implementation of the Whitham procedure in the analysis of oscillatory behavior of solutions to the KdV equations). The latter will also be applied to the asymptotic description of solutions inside the oscillatory zone.
Furthermore we plan to study the possibility of extending the Main Conjecture to the critical behavior of solutions to the Hamiltonian perturbations of more general first order quasilinear systems of elliptic type. Last but not least, it would be of interest to study the distribution of poles of the \emph{tritronqu\'ee} solution in the sector $|\arg \zeta| >\frac{4\pi}{5}$ and to compare these poles with the peaks of solutions to NLS inside the oscillatory zone. The elliptic asymptotics obtained by Kitaev \cite{kita} might be useful for studying these poles for large $|\zeta|$.

In this paper we did not study the behaviour of solutions to NLS near the boundary $u=0$. Such a study is postponed for a subsequent publication.



\newpage

\begin{thebibliography}{99}

\bibitem{ag} G.P.Agrawal, {\it Nonlinear Fiber Optics}. Academic Press, San Diego, 2006, 4th edition.

\bibitem{al} S.Alinhac, {\it Blowup for Nonlinear Hyperbolic Equations.} Progress in Nonlinear Differential Equations and their Applications, {\bf 17}. Birkh\"auser Boston, Inc., Boston, MA, 1995.

\bibitem{ar} V.I.Arnold,V.V.Goryunov, O.V.Lyashko,V.A.Vasil'ev,  {\it Singularity Theory. I.}  Dynamical systems. VI, Encyclopaedia Math. Sci. {\bf 6}, Springer, Berlin, 1993. 

\bibitem{bo} P.Boutroux, Recherches sur les transcendants de M. Painlev\'e et l'\'etude asymptotique des \'equations diff\'erentielles du second ordre. {\it Ann. \'Ecole Norm} {\bf 30} (1913) 265 - 375.

\bibitem{bronski} J.C.Bronski, J.N.Kutz, 
Numerical simulation of the semiclassical limit of the focusing nonlinear Schr\"odinger equation.
{\it Phys. Lett.} {\bf A 254} (2002) 325 - 336.

\bibitem{venak} R.Buckingham, S.Venakides, Long-time asymptotics of the nonlinear Schr\"odinger equation shock problem. {\it Comm. Pure Appl. Math.}, Published Online 12.03.2007.

%
\bibitem{ca} R.Carles, WKB analysis for the nonlinear Schr\"odinger equation and instability results. ArXiv:math.AP/0702318.


\bibitem{ce} H.D.Ceniceros, F.-R.Tian, A numerical study of the semi-classical limit of the focusing nonlinear Schr\"odinger equation. 
{\it Phys. Lett.} {\bf A 306} (2002) 25--34. 

\bibitem{cl} T.Claeys, M.Vanlessen, The existence of a real pole-free solution of the fourth order analogue of the Painlev\'e I equation. ArXiv:math-ph/0604046.

\bibitem{co} O.Costin, Correlation between pole location and asymptotic behavior for Painlev\'e I solutions. {\it Comm. Pure Appl. Math.} {\bf 52} (1999) 461--478. 

\bibitem{ch} M.C.Cross and P.C.Hohenberg, Pattern formation outside of equilibrium. {\it Rev. Mod. Phys.} {\bf 65} (1993) 851-1112.

\bibitem{du1} B.Dubrovin, S.-Q.Liu, Y.Zhang, On Hamiltonian perturbations of hyperbolic systems of conservation laws I: quasitriviality of bihamiltonian perturbations. {\it Comm. Pure Appl. Math.} {\bf 59} (2006) 559-615.

\bibitem{du2} B.Dubrovin, On Hamiltonian perturbations of hyperbolic systems of conservation laws, II: universality of critical behaviour, 
{\it Comm. Math. Phys.} {\bf 267} (2006) 117 - 139.

\bibitem{duits} M.Duits, A.Kuijlaars, Painlev\'e I asymptotics for orthogonal polynomials with respect to a varying quartic weight.
ArXiv:math/0605201.

\bibitem{fokas} A.S.Fokas, S.Tanveer, A Hele - Shaw problem and the second Painlev\'e transcendent. {\it Math. Proc. Camb. Phil. Soc.} {\bf 124} (1998) 169 - 191.

\bibitem{fo}  M.G.Forest, J.E.Lee, Geometry and modulation theory for the periodic nonlinear Schr\"odinger equation. In: {\it Oscillation Theory, Computation, and Methods of Compensated Compactness}
(Minneapolis, Minn., 1985), 35-69. The IMA Volumes in Mathematics and Its Applications, 2. 
Springer, New York, 1986. 

%
\bibitem{gwp1} J.~Ginibre, G.Velo, On a class of nonlinear Schr\"odinger equations. I. The Cauchy problem, general case. {\it J. Funct. Anal.} {\bf 32} (1979) 1-32.
    
\bibitem{gradrhy}     
I. S.Gradshteyn, I. M. Ryzhik, {\it Table of Integrals, Series, and Products.} Translated from the Russian. Sixth edition. Translation edited and with a preface by Alan Jeffrey and Daniel Zwillinger. Academic Press, Inc., San Diego, CA, 2000.

\bibitem{gk1} T.Grava, C.Klein, Numerical solution of the small dispersion limit of Korteweg de Vries and Whitham equations. ArXiv:math-ph0511011, to appear in {\it Comm. Pure Appl. Math.}, 2007.

\bibitem{gk} T.Grava, C.Klein, Numerical study of a multiscale expansion of KdV and Camassa-Holm equation. ArXiv:math-ph/0702038.

\bibitem{gr} E.Grenier, Semiclassical limit of the nonlinear Schr\"odinger equation in small time. {\it Proc. Amer. Math. Soc. } {\bf 126}  (1998) 523--530.

\bibitem{in} E.L.Ince, {\it Ordinary Differential Equations}. Dover Publications, New York, 1944.

\bibitem{jin} S.Jin, C.D.Levermore, D.W.McLaughlin, The behavior of solutions of the NLS equation in the semiclassical limit. {\it Singular Limits of Dispersive Waves} (Lyon, 1991), 235--255, NATO Adv. Sci. Inst. Ser. B Phys., {\bf 320}, Plenum, New York, 1994.


\bibitem{jk} N.Joshi, A.Kitaev, On Boutroux's tritronqu\'ee solutions of the first Painlev\'e equation. {\it Stud. Appl. Math.} {\bf 107} (2001) 253--291.

\bibitem{kam1} S.Kamvissis, Long time behavior for the focusing nonlinear Schr\"odinger equation with real spectral singularities.  {\it Comm. Math. Phys.} {\bf 180}  (1996) 325--341.

\bibitem{kam2} S. Kamvissis, K.D.T.-R.McLaughlin, P.D.Miller,  {\it Semiclassical Soliton Ensembles for the Focusing Nonlinear Schr\"odinger Equation.} Annals of Mathematics Studies, 154. Princeton University Press, Princeton, NJ, 2003.

\bibitem{ka} A.Kapaev, Quasi-linear Stokes phenomenon for the Painlev\'e first equation. {\it J. Phys.  A: Math. Gen.} {\bf 37} (2004) 11149--11167.

\bibitem{kita} A.Kitaev, The isomonodromy technique and the elliptic asymptotics of the first Painlev\'e transcendent.  {\it Algebra i Analiz} {\bf 5} (1993), no. 3, 179--211; translation in 
{\it St. Petersburg Math. J.} {\bf 5} (1994), no. 3, 577--605.

\bibitem{numart1d} C. Klein, Fourth order time-stepping for low dispersion Korteweg - de Vries and nonlinear Schr\"odinger equation (2006),
http://www.mis.mpg.de/preprints/2006/prepr$2006\_{}$133.html

\bibitem{krasny} R.Krasny, A study of singularity formation in a vortex sheet by the point-vortex approximation.  {\it J. Fluid Mech.}  {\bf 167}  (1986) 65--93.

\bibitem{optim}J.~C.~Lagarias, J.~A.~Reeds, M.~H.~Wright, 
 and P.~E.~Wright, Convergence properties of the Nelder-Mead 
 simplex method in low dimensions. {\it SIAM Journal of Optimization} 
{\bf 9} (1988) 112-147. 

\bibitem{miller} G.D.Lyng, P.D.Miller, The $N$-soliton of the focusing nonlinear Schr\"odinger equation for $N$ large.
{\it Comm. Pure Appl. Math.} {\bf 60} (2007) 951-1026.

\bibitem{me} G.M\'etivier, Remarks on the well-posedness of the nonlinear Cauchy problem. ArXiv:math.AP/0611441.

\bibitem{mi} P.D.Miller, S.Kamvissis, On the semiclassical limit of the focusing nonlinear Schr\"odinger equation. {\it  Phys. Lett.} {\bf A  247}  (1998) 75--86.

\bibitem{new} A.C.Newell, {\it Solitons in Mathematics and Physics}. CBMS-NSF Regional Conference Series in Applied Mathematics, {\bf 48}. SIAM, Philadelphia, PA, 1985. 

\bibitem{novik} S.P.Novikov, S.V.Manakov, L.P.Pitaevski\u\i, V.E.Zakharov, {\it Theory of Solitons. The Inverse Scattering Method.} Translated from the Russian. Contemporary Soviet Mathematics. Consultants Bureau [Plenum], New York, 1984.

\bibitem{saya}\textsc{J.~Satsuma and N.~Yajima}, \emph{Initial value 
problems of one-dimensional self-modulation of nonlinear waves in 
dispersive 
media}, Supp. Prog. Theo. Phys. 55 (1974), pp.~284Ð-306. 

\bibitem{sh} A.B.Shabat,  One-dimensional perturbations of a differential operator, and the inverse scattering 
problem. In: {\it Problems in Mechanics and Mathematical Physics}, 279--296. Nauka, Moscow, 1976. 

\bibitem{bvp4c}  L.~F.~Shampine,  M.~W.~Reichelt and J.~Kierzenka, 
\textit{Solving Boundary Value Problems for Ordinary Differential 
Equations in MATLAB with bvp4c}, available at 
http://www.mathworks.com/bvp\_tutorial

\bibitem{si} P.Sikivie, The caustic ring singularity. {\it Phys. Rev.} {\bf D60} (1999) 063501.

\bibitem{sl} M.Slemrod, Monotone increasing solutions of the Painlev\'e 1 equation $y''=y\sp 2+x$ and their role in the stability of the plasma-sheath transition. {\it European J. Appl. Math.} {\bf 13} (2002) 663--680.

\bibitem{th} R.Thom, {\it Structural Stability and Morphogenesis: An Outline of a General Theory of Models.} Reading, MA: Addison-Wesley, 1989.

\bibitem{to1} A.Tovbis, S.Venakides, X.Zhou, On semiclassical (zero dispersion limit) solutions of the focusing nonlinear Schr\"odinger equation. {\it Comm. Pure Appl. Math.} {\bf  57}  (2004) 877--985.

\bibitem{to2} A.Tovbis, S.Venakides, X.Zhou,  On the long-time limit of semiclassical (zero dispersion limit) solutions of the focusing nonlinear Sch\"odinger equation: pure radiation case.  {\it Comm. Pure Appl. Math.}  {\bf 59}  (2006) 1379--1432.

\bibitem{trefethen1}  L.~N.~Trefethen,
   \emph{ Spectral Methods in MATLAB}, SIAM, Philadelphia, PA, 2000.
   
\bibitem{gwp2} Y.Tsutsumi, $L^2$-solutions for nonlinear Schr\"odinger equations and nonlinear groups. {\it Funkcial. Ekvac.} {\bf 30} (1987) 115--125.

\bibitem{whi} G.B.Whitham, {\it Linear and Nonlinear Waves}. Wiley-Intersci. 1974.

\bibitem{za} V.E.Zakharov, A.B.Shabat, A. B. Exact theory of two-dimensional self-focusing and one-dimensional self-modulation of waves in nonlinear media. {\it Soviet Physics JETP} {\bf 34} (1972), 
no. 1, 62-69.; translated from {\it \v{Z}. Eksper. Teoret. Fiz.}  (1971), no. 1, 118-134.

\end{thebibliography}
\end{document}